\documentclass[12pt,onecolumn,oneside,leqno,a4paper]{article}
\pdfoutput=1
\usepackage{german}
\usepackage{exscale}
\usepackage{amsmath,amssymb,verbatim,eucal} 
\usepackage{amscd}
\usepackage{graphicx}
\usepackage{rotating}
\usepackage{amsxtra,amsthm}
\usepackage{color}
\usepackage{makeidx}
\usepackage{xspace} 
\usepackage{framed} 
\usepackage[pdftex]{hyperref} 
\usepackage{xr} 
\usepackage{mathtools} 

\newcommand{\BT}{\begin{theorem}}
\newcommand{\ET}{\end{theorem}}
\newcommand{\TB}[1]{\textbf{#1}}
\newcommand{\TI}[1]{\textit{#1}}

\newcommand\QED{\hspace*{\fill} $ \square $}

\newenvironment{DES}[1]
    {\begin{list}{}
    {\settowidth{\labelwidth}{#1}
	\setlength{\labelsep}{1em} 
	\setlength{\leftmargin}{\labelwidth}
	\addtolength{\leftmargin}{\labelsep}
	\setlength{\parsep}{\parskip}
	\setlength{\itemsep}{0.2\parsep plus0.02\parsep minus0.02\parsep}
	\setlength{\topsep}{0.2\parsep plus0.02\parsep minus0.02\parsep}
	\setlength{\partopsep}{0.6\parsep plus0.02\parsep minus0.05\parsep}
	}}
	{\end{list}}
	
\newenvironment{DES1}[1]
	{\begin{list}{}
	{\settowidth{\labelwidth}{#1}
	\setlength{\labelsep}{1em} 
	\setlength{\leftmargin}{0em}
	\setlength{\itemindent}{\labelwidth}
	\addtolength{\itemindent}{\labelsep}
	\setlength{\parsep}{\parskip}
	\setlength{\itemsep}{0.2\parsep plus0.02\parsep minus0.02\parsep}
	\setlength{\topsep}{0.2\parsep plus0.02\parsep minus0.02\parsep}
	\setlength{\partopsep}{0.6\parsep plus0.02\parsep minus0.05\parsep}
	}}
	{\end{list}}

\def\R{\mathbb{R}}

\def\Rplus{\mathbb{R}^+}

\def\B{\mathbb{B}}

\def\SS{\mathbb{S}}

\def\diam{\text{\textup{diam}}}

\def\const{\text{\textup{const.}}}

\newcommand\Forall{\;\forall\;}

\newcommand\Min{\min\,}

\newcommand\epsi{\varepsilon}
\newcommand\vhi{\varphi}

\renewcommand\Dot[2]{\left\langle #1,#2 \right\rangle}

\def\ds{\displaystyle}

\renewcommand{\Tilde}[1]{\widetilde{\,#1\,}}

\newcommand{\Tildeu}[1]{\!\widetilde{\,#1\,}{}}

\newcommand{\To}{\longrightarrow}

\newcommand{\Ast}{{\ds \ast}}
\newcommand{\Circ}{{\ds \circ}} 

\newcommand{\upto}{\uparrow }

\newcommand{\BM}{\begin{pmatrix}}
\newcommand{\EM}{\end{pmatrix}}

\newtheoremstyle{tld}
{}
{}
{\itshape}
{}
{\bfseries}
{.}
{2em}
{}

\newtheoremstyle{bb}
{}
{}
{}
{}
{\bfseries}
{.}
{2em}
{}

\newtheoremstyle{bbs}
{}
{}
{}
{}
{\bfseries}
{.}
{\newline}
{}

\swapnumbers

\usepackage{xpatch}
\makeatletter
\xpatchcmd{\@thm}{\fontseries\mddefault\upshape}{}{}{}
\makeatother

\theoremstyle{tld}
\newtheorem{thm}{\negthickspace. Theorem}[section] 
\newtheorem{thmz}[thm]{\negthickspace. Supplement}
\newtheorem{lem}[thm]{\negthickspace. Lemma}
\newtheorem{cor}[thm]{\negthickspace. Corollary}
\newtheorem{defi}[thm]{\negthickspace. Definition}

\theoremstyle{bb}
\newtheorem{ex}[thm]{\negthickspace. Example}
\newtheorem{rem}[thm]{\negthickspace. Remark}
\newtheorem{noth}{}[section]

\theoremstyle{bbs}
\newtheorem{exs}[thm]{\negthickspace. Examples}
\newtheorem{rems}[thm]{\negthickspace. Remarks}

\renewcommand{\BT}{\begin{thm}}
\renewcommand{\ET}{\end{thm}}

\newcommand{\BTZ}{\begin{thmz}}
\newcommand{\ETZ}{\end{thmz}}

\newcommand{\BL}{\begin{lem}}
\newcommand{\EL}{\end{lem}}

\newcommand{\BC}{\begin{cor}}
\newcommand{\EC}{\end{cor}}

\newcommand{\BD}{\begin{defi}}
\newcommand{\ED}{\end{defi}}

\newcommand{\BX}{\begin{ex}}
\newcommand{\EX}{\end{ex}}

\newcommand{\BXS}{\begin{exs}}
\newcommand{\EXS}{\end{exs}}

\newcommand{\BR}{\begin{rem}}
\newcommand{\ER}{\end{rem}}

\newcommand{\BRS}{\begin{rems}}
\newcommand{\ERS}{\end{rems}}

\newcommand{\BN}{\begin{noth}}
\newcommand{\EN}{\end{noth}}

\newcommand{\BE}{\begin{equation}}
\newcommand{\EE}{\end{equation}}

\newcommand{\abs}[1]{\left|#1\right|}

\newcommand{\vac}{\varnothing}

\def\l{\ell}
\def\ls{\mspace{1.0mu}}

\numberwithin{equation}{section}

\newcommand{\bild}[1]{\includegraphics{#1}}

\newcommand{\cross}{\times}

\begin{document}

\begin{center}
\TB{\LARGE On a Minimax Problem for Ovals}

\TB{\large Rolf Walter}
\end{center}

\setcounter{section}{0}


\section{\hspace{-1em}.
Introduction}
\label{intro_1}
\markright{\ref{intro_1}. Introduction}


\thispagestyle{empty}
Let $ (X,d) $ be a bounded metric space. A
fundamental invariant of $ X $ is the \TI{diameter}\\[-3ex]
\BE
\diam(X) := \text{sup}_{p,q \in X} \; d(p,q).
\label{eq_1}
\EE
Besides $ \diam(X) $, one can form the quantity
\BE
\delta(X) := 
\text{inf\rule[-0.5ex]{0em}{1ex}}_{\,p\in X}\;
\text{sup}_{q \in X} \; d(p,q),
\label{eq_2}
\EE
which, in general, is distinct from $ \diam(X) $. 

$ \delta(X) $ is the smallest radius of balls covering $ X $ with centers in $ X $. Despite this obvious meaning, the invariant 
$ \delta(X) $ seems rather untouched in the geometric literature.
In the area of approximation theory, $ \delta(X) $ occurs as the relative Chebyshev radius of $ X $ w.r.t. $ X $ itself, i.e. 
$ \delta(X) = r(X,X) $ in the terminology of \TI{Amir/Ziegler} [1980]. The size of $ \delta(X) $ might have an impact on
the shape of the space $ X $. In particular, for convex hypersurfaces in $ \R^{n} $, there could exist relations to the many other geometric invariants of such objects. 

In this paper we
discuss, for plane convex curves $ X $, an isoperimetric
type inequality, relating $ \delta(X) $ to the perimeter of
the curve. Due to the minimax character of $ \delta(X) $, 
its handling resists the usual principles for extremal 
questions. It will be shown that the smooth case of the 
inequality can be reduced to the polygonal case 
via approximation. However, for polygons, there is the 
additional problem of the high dimensionality of the set 
of vertices. So, in
general, we only can offer conjectures. Definite
solutions are possible for restricted classes of curves.
Even for short polygons, there already arises a sort 
of `magic kites' which show that squares are definitely not 
the extremal figures.

In the general context, it is easy to see that the function of farthest distance 
\BE
\mu: X \To \R, \quad 
\mu(p) := \text{sup}_{q \in X} \; d(p,q)
\label{eq_2a}
\EE
is Lipschitz-continuous with Lipschitz-constant $ 1 $. In fact, from the triangle inequalities for points $ p_{1},p_{2},q \in X $:
\BE
d(p_{1},q)-d(p_{1},p_{2}) \leq
d(p_{2},q) 
\leq d(p_{1},q)+d(p_{1},p_{2})
\label{eq_3b}
\EE
one derives, by taking suprema w.r.t. $ q $
\BE
\mu(p_{1})-d(p_{1},p_{2}) \leq
\mu(p_{2}) 
\leq\mu(p_{1})+d(p_{1},p_{2}),
\label{eq_3c}
\EE
hence
\BE
|\mu(p_{1})-\mu(p_{2})| \leq d(p_{1},p_{2}).
\label{eq_3d}
\EE

If $ X $ is compact then the suprema in \eqref{eq_1} and
\eqref{eq_2} may be replaced equivalently by the maxima,
because $ d(p,q) $ is continuously dependent on the argument
pair. By the continuity of $ \mu $, the same is true for the
infimum in \eqref{eq_2}. More precisely we can state:

\BL
\label{lem_1} 
Let $ X $ be a compact metric space $ X $. Then we have
\BE
\delta(X) = 
\textup{min}_{\,p\in X}\;\textup{max}_{\,q \in X} \;d(p,q).
\label{eq_3}
\EE
In particular, there are $ p_{0},q_{0} \in X $ with 
\BE
\delta(X) = d(p_{0},q_{0}) \geq d(p_{0},q) \quad 
\text{for all $ q \in X $}.
\label{eq_3a}
\EE
Conversely, if $ p_{0},q_{0} \in X $ are such that 
\eqref{eq_3a} holds, then the function 
$ \mu(p) := \sup_{q \in X} d(p,q) $ 
$ = \max_{q \in X} d(p,q) $ 
attains its minimum at $ p_{0} $, and the function 
$ q \mapsto d(p_{0},q) $ attains its maximum at $ q_{0} $.
\EL

\TI{Proof.} 
The function $ \mu $ assumes its minimum at some point 
$ p_{0} \in X $, and then the function
$ r \mapsto d(p_{0},r) $ assumes its maximum at some point 
$ q_{0} \in X $. Thus
\BE
\text{max}_{\,r \in X} \; d(p,r) \geq 
\text{max}_{\,r \in X} \; d(p_{0},r) = 
d(p_{0},q_{0}) \geq 
d(p_{0},q) \quad \text{for all $ q \in X $}.
\label{eq_3e}
\EE
This implies by switching to the infimum w.r.t. $ p $
\BE
\delta(X) = 
\textup{min}_{\,p \in X}\;\textup{max}_{\,r \in X} \; d(p,r) = 
d(p_{0},q_{0}) \geq 
d(p_{0},q) \quad \text{for all $ q \in X $}.
\label{eq_3f}
\EE
This proves \eqref{eq_3} and \eqref{eq_3a} for this 
specific pair $ (p_{0},q_{0}) $.

Conversely, if $ p_{0},q_{0} $ is \TI{any} pair in $ X $, 
satisfying \eqref{eq_3a}, then 
$ \min_{p \in X} \mu(p) =  d(p_{0},q_{0}) \geq d(p_{0},q) $ 
for all $ q \in X $. This implies 
$ \mu(p_{0}) = \max_{q \in X} d(p_{0},q) = d(p_{0},q_{0}) $, 
hence $ \min_{p \in X} \mu(p) = \mu(p_{0}) $. The second
relation in \eqref{eq_3a} shows the remainder of the
assertion. \QED

We call a pair of points $ p_{0},q_{0} \in X $ 
\TI{distinguished} if it fulfills \eqref{eq_3a}. 
Equivalent to this is that the function $ \mu $ has the
minimum point $ p_{0} $ and the function $ q \mapsto
d(p_{0},q) $ has the maximum point $ q_{0} $. So 
$ p_{0},q_{0} $ can well be viewed as points where the
minimax-invariant $ \delta(X) $ attains its value. Observe
that the roles of $ p_{0},q_{0} $ are not interchangeable.
It is possible that there are more than one pair of
distinguished points, e.g. if $ X $ has isometric
symmetries. In the terminology of approximation theory, 
the first point $ p_0 $ of a distinguished pair $ p_{0},q_{0} $ is the same as a relative Chebyshev center.

Now let $ X $ be a simply closed convex curve $ \Gamma $ in 
$ \R^{2} $ equipped with the metric $ d $, restricted from the
standard Euclidean metric in $ \R^{2} $. Assume that the 
domain $ \Omega $ enclosed by $ \Gamma $ is non-void. 
Denote by $ C := \Omega \cup \Gamma $ the corresponding
convex body and by $ L(\Gamma) $ its perimeter.

A distinguished pair of points $ p_{0},q_{0} $ at which 
$ \delta $ assumes its value according to \eqref{eq_3a} determines the related chord as the segment from  $ p_{0} $ to $ q_{0} $. Thereby, $ q_{0} $ is the maximal point of the function 
$ q \mapsto d(p_{0},q) $. Such a chord (which of course is contained in $ C $) will be called \TI{distinguished}. Thus, each closed convex curve determines certain peculiar chords which, without doubt, have a geometric meaning. However, at this time, only few is known on this. E.g. for a $ C^{2} $-curve, it is open how the distinguished chords can be determined. For convex polygons we shall detect some details of distinguished chords in section \ref{alg_4a}. For instance, the second endpoint of a distinguished chord must be a vertex, and there exists a finite algorithm for the finding of all distinguished chords.

Since each circle of radius $ \delta(\Gamma) $ around a relative 
Chebyshev center $ p_0 $ encloses $ \Gamma $, one has
$$
L(\Gamma) \leq 2\pi \delta(\Gamma),
$$
by the monotony of the boundary volume for convex sets (see 
\TI{Bangert} [1981] for a rather general version of the monotony). In the opposite direction we state the

\TB{Conjecture 1:} 
\TI{
For any closed convex curve $ \Gamma $ in the plane $ \R^{2} $, the inequality
\BE
L(\Gamma) \geq \pi \cdot \delta(\Gamma)
\label{eq_3.1}
\EE 
holds true.
}

This conjecture is supported by the results of this paper and by numerous experiments on further classes of curves. 

Returning to the definition of 
$ \delta(\Gamma) $, one can state the conjecture \eqref{eq_3.1} as follows: \TI{Let $ \delta_{0} $ be a positive constant such that, for any point $ p \in \Gamma $, there is a point 
$ q \in \Gamma $ with $ d(p,q) \geq \delta_{0} $. Then 
$ L(\Gamma) \geq \pi \cdot \delta_{0} $.}
If this were true in general it would answer a question recently raised in the field of Reverse Engineering (\TI{Wortmann} [2004]).


\section{\hspace{-1em}. Curves nearby Circles} 
\label{cnc}
\markright{\ref{cnc}. Curves nearby Circles}


By this we will understand that all the curvature centers 
fall inside the curve:

\BT
\label{thm_2}
Let the closed convex curve $ \Gamma $ in $ \R^{2} $ be  
regular and of class $ C^{2} $. Assume that all curvature 
centers of $ \Gamma $ lie in the interior $ \Omega $. Then 
we have
\BE
L(\Gamma) \geq \pi \cdot \delta(\Gamma).
\label{eq_7}
\EE
Equality in \eqref{eq_7} holds if and only if $ \Gamma $ is 
of constant breadth.
\ET

For any $ p \in \Gamma $ there is at least one point 
$ q_{1} \in \Gamma $ where the function 
$ q \mapsto d(p,q) $ reaches its maximum. The closed segment 
$ [pq_{1}] $ is orthogonal to the tangent of $ \Gamma $ at 
$ q_{1} $. However, it is not clear whether $ q_{1} $ is
unique or how it depends on $ p $. We therefore consider an
eventual inverse of this map, defined as follows: For any 
$ q \in \Gamma $ let $ \nu(q) $ be the inward normal of 
$ \Gamma $ at $ q $ (counted as a half-ray starting at 
$ q $). It cuts $ \Gamma $ in a unique second point 
$ \Phi(q) $, and the open segment $ ]q,\Phi(q)[ $ is 
contained in the interior $ \Omega $. The map
\BE
\Phi: \Gamma \To \Gamma
\label{eq_2.1}
\EE
will play a decisive role in our arguments.

The orientation of $ \Gamma $ shall always be such that the 
interior $ \Omega $ lies to the left. At any point 
$ p \in \Gamma $ we then have a unique inner unit normal  
vector $ N(p) $ of $ \Gamma $, and for any positively 
oriented tangent vector $ w $ of $ \Gamma $ at $ p $, the 
pair $ (w,N(p)) $ will be positively oriented w.r.t 
$ \R^{2} $.

Now, let the segment $ [p,q] $ be any \TI{chord} of 
$ \Gamma $, i.e. $ p, q $ are distinct points on $ \Gamma $,  
and $ ]p,q[ $ belongs to the interior $ \Omega $. The 
straight line $ g $ spanned by $ p,q $ decomposes $ \Gamma $ 
into two arcs which lie on different sides of $ g $. For
each of these arcs, $ p $ and $ q $ are the points where the
arc enters or leaves (transversally) the half plane of $ g $
which contains the arc. So, if $ u $, resp. $ v $ are
tangent vectors of $ \Gamma $ at $ p $, resp. $ q $ then 
$ u,v $ are equally oriented w.r.t. $ \Gamma $ if and only if
they point into different sides of $ g $. 

\bigskip
\BL
\label{lem_2}
Under the assumptions of Theorem \ref{thm_2}, the map 
$ \Phi $ in \eqref{eq_2.1} is an orientation-preserving 
$ C^{1} $-diffeomorphism of $ \Gamma $ onto itself.
\EL

\medskip
\TI{Proof.} 
Let $ \gamma: [0,\l\ls] \to \R^{2} $ be a positive 
parametrization of $ \Gamma $ by arclength $ s $ such that 
$ \gamma $ is injective on $ [0,\l[ $ and 
$ \gamma(0) = \gamma(\l) $. (If necessary, we may extend 
$ \gamma $ periodically to $ \R $ with simple period 
$ \l $.) For the parametrizations of the unit tangent and
normal vectors $ T $ and $ N $ we have the structure
equations of Darboux
\BE
\gamma{\,}' = T, \qquad
T' = \kappa N, \qquad
N' = - \kappa T,
\label{eq_2.2}
\EE
$ \kappa $ being the curvature function. The curvature
centers are given by
\BE
z = \gamma+\frac{N}{\kappa},
\label{eq_2.3}
\EE
where $ \kappa $ is always positive because the curvature 
centers are in $ \Omega $ (thus finite): $ \Gamma $ is 
strictly convex. For the $ \Phi $-images of the curve points 
we have a representation of the form
\BE
\Phi(\gamma(s)) = \gamma(s)+D(s)N(s) = \gamma(\vhi(s)),
\label{eq_2.4}
\EE
where $ D(s) $ is unique, namely equal to the distance 
$ d(\gamma(s),\Phi(\gamma(s)) $, and $ \vhi(s) $ is unique 
modulo $ \l $. 
In order to show the differentiability of 
these quantities we consider the $ C^{1} $-function
\BE
F(s,\xi,\eta) := \gamma(s)+\xi N(s) - \gamma(\eta)
\label{eq_2.6}
\EE
and apply the implicit function theorem to the equation 
$ F(s,\xi,\eta) = 0 $. The partials of $ F $ w.r.t. $ \xi $
and $ \eta $ are $ N(s) $ and $ -\gamma{\,}'(\eta) $ which
for $ \eta = \vhi(s) $ are linearly independent by the
transversality mentioned above. Thus $ D $ is of class 
$ C^{1} $, and the same is true for $ \Phi $ by
\eqref{eq_2.4}. Then $ \vhi $ can be chosen continuously,
satisfying $ \Phi \circ \gamma = \gamma \circ \vhi $, so in
fact being of class $ C^{1} $. (Here it is necessary to view
all functions on the whole of $ \R $.)

The assumption on the curvature centers is 
expressed by 
\BE
D(s) > \frac{1}{\kappa(s)} \quad \text{for all $ s $}.
\label{eq_2.5}
\EE
Differentiating the right equation in \eqref{eq_2.4} and 
inserting according to \eqref{eq_2.2} we obtain
\BE
(1-D\kappa)T+D'N = (\gamma{\,}' \circ \vhi)\cdot \vhi'
\label{eq_2.7}
\EE
and by scalar multiplication with $ \gamma{\,}' $
\BE
(1-D\kappa)\Dot{\gamma{\,}'}{\gamma{\,}'} = 
\Dot{\gamma{\,}' \circ \vhi}{\gamma{\,}'} \cdot \vhi'.
\label{eq_2.8}
\EE
For each $ s $, $ \gamma{\,}'(s) $ is a normal vector of the 
chord $ [\gamma(s),\Phi(\gamma(s))] $, so the two scalar
products in \eqref{eq_2.8} have different signs. On the
other hand, $ 1-D\kappa $ is always negative by
\eqref{eq_2.5}; thus, by \eqref{eq_2.8}, $ \vhi' $ is always
positive, hence $ \vhi $ orientation-preserving.

This proves that $ \Phi $ is an immersion of the compact 
manifold $ \Gamma $ into $ \Gamma $, hence a covering 
projection (\TI{Kobayashi/Nomizu} [1963], Cor. 4.7), say with 
$ m \geq 1 $ leaves. $ \Gamma $ is homeomorphic to the 
sphere $ \SS^{1} $. Since $ \Phi $ is orientation-preserving,
its mapping degree also equals $ m $. If we had $ m > 1 $
then, by the Hopf/Lefschetz fixpoint theorem, there would
exist a fixpoint of $ \Phi $ what is not possible. So 
$ \Phi $ is injective and has all properties, as asserted.
\QED

\bigskip
\TI{Proof of Theorem \ref{thm_2}}

\medskip
We continue in the setting of the last proof. The mapping
$ \Phi \circ \gamma = \gamma \circ \vhi $ is now called 
$ \Tilde{\gamma} $ and viewed on $ [0,\l\ls] $. As $ \gamma $ 
itself, $ \Tilde{\gamma} $ is a positive parametrization of 
$ \Gamma $, injective on $ [0,\l\ls [ $ with 
$ \Tilde{\gamma}(0) = \Tilde{\gamma}(\l) $. By 
\eqref{eq_2.4}, $ \Tilde{\gamma} $ and its derivative 
are given by
\BE
\Tilde{\gamma} = \gamma+DN, \qquad
\Tildeu{\gamma}' = (1-D\kappa)T+D'N.
\label{eq_2.9}
\EE
This implies for the length of $ \Gamma $:
\BE
\l 
= \int_{0}^{\l} |\Tildeu{\gamma}'| 
= \int_{0}^{\l} \sqrt{(D\kappa-1)^{2}+D'{}^{2}} 
\geq \int_{0}^{\l} (D\kappa-1) 
= \int_{0}^{\l} D\kappa - \l.
\label{eq_2.10}
\EE
For each $ \Tilde{\gamma}(s) =: p $, the $ \gamma(s) $ is
the only point $ q \in \Gamma $ where the distance 
$ d(p,q) $ is maximal. Namely, if $ q_{1} \in \Gamma $ were
a point $ \neq \gamma(s) $ with maximal distance from $ p $
then the segment $ [p,q_{1}] $ were orthogonal to $ \Gamma $
at $ q_{1} $, hence $ \Phi(q_{1}) = p $, in contradiction to
the injectivity of $ \Phi $. This implies 
\BE
D(s) = 
\text{max}_{\,q \in \Gamma} \; d(\Tilde{\gamma}(s),q)
\geq \delta(\Gamma).
\label{eq_2.11}
\EE
Thus from \eqref{eq_2.10} we infer, using finally the 
theorem on turning tangents:
\BE
2 \l 
\geq \delta(\Gamma) \int_{0}^{\l} \kappa
= \delta(\Gamma) \cdot 2\pi.
\label{eq_2.12}
\EE
This is the assertion \eqref{eq_7}. 

\medskip
If we have equality in \eqref{eq_7} then we deduce from  
\eqref{eq_2.10} to \eqref{eq_2.12} that 
$ D = \delta(\Gamma) $. For constant $ D $, \eqref{eq_2.9} 
implies $ \Tildeu{\gamma}' = (1-D\kappa)T $, so the tangents
of $ \Gamma $ at $ \gamma(s) $ and $ \Tilde{\gamma}(s) $ are
parallel and 
$ d(\gamma(s),\Tilde{\gamma}(s)) = D = \const $. This shows 
that $ \Gamma $ has constant breadth $ \delta(\Gamma) $.

On the other hand, if a closed regular and convex 
$ C^{2} $-curve $ \Gamma $ has constant breadth $ a $ then,
by well-known properties of such ovals, the curvature
centers are situated in the interior and 
$ \delta(\Gamma) = a = L(\Gamma)/\pi $; see e.g. 
\TI{Bonnesen/Fenchel} [1934], Sect. 63, p. 128. \QED

\newpage

\section{\hspace{-1em}. Approximating \texorpdfstring{
$ \pmb{\delta}\pmb{(}\mathbf{\Gamma}\pmb{)} $}{deltaX}}
\label{appr_3}
\markright{\ref{appr_3}. Approximating $ \delta(\Gamma) $}


An important question is whether the minimax-invariant 
$ \delta(\Gamma) $ could be calculated via polygonal
approximations. Theorem \ref{lem_3.4} below provides
a result in this direction. It contains an error estimate 
which may be useful for numerical purposes. 
But above all, eventual improvements of the inequality \eqref{eq_7} are thus reduced to the polygonal case. 

First, we observe two simple facts:

\medskip
\BL
\label{lem_3.2}
~\\
Let $ X, Y $ be bounded subsets of a metric space $ Z $ 
(with the restricted metrics).
\begin{DES}{(ii)}
\item[(i)]
Let  $ \epsi $ be a positive
constant. If there exists a surjective map $ \vhi: X \to Y $
such that $ d(p,\vhi(p)) \leq \epsi $ for all $ p \in X $
then $ \delta(X) \leq \delta(Y)+2\epsi $.
\item[(ii)]
If there is a surjective map $ \vhi: X \to Y $ satisfying
$ d(p,q) \leq d(\vhi(p),\vhi(q)) $ for all $ p,q \in X $, 
then $ \delta(X) \leq \delta(Y) $.
\end{DES}
\EL

\TI{Proof.} \\
\TI{For (i):} 
For any points $ p,q \in X $ we have
\BE
d(p,q) \leq 
d(p,\vhi(p))+d(\vhi(p),\vhi(q))+d(q,\vhi(q)) \leq
d(\vhi(p),\vhi(q))+2\epsi,
\label{eq_3.4}
\EE
thus, taking suprema w.r.t. $ q $ for fixed $ p $,
\BE
\text{sup}_{q \in X} \; d(p,q) \leq
\text{sup}_{q \in X} \; d(\vhi(p),\vhi(q))+2\epsi =
\text{sup}_{q_{1} \in Y} \; d(\vhi(p),q_{1})+2\epsi,
\label{eq_3.5}
\EE
since $ \vhi $ is surjective. For the same reason, taking 
infima w.r.t. $ p $, this yields the assertion.

\TI{For (ii):} The assertion is immediately clear from the 
hypothesis, taking first suprema w.r.t. $ q $ and then 
infima w.r.t. $ p $, observing that $ \vhi(q) $ and 
$ \vhi(p) $ run through all of $ Y $. \QED

\bigskip
In case of an oriented closed convex and regular 
$ C^{2} $-curve $ \Gamma $, consider a chord $ [p,q] $ and
orient the line $ g $ connecting $ p,q $ in such a way that
its right half plane $ H_{g}^{+} $ is entered by the curve at
$ p $ and left at $ q $. The convex body 
$ C := \Omega \cup \gamma $ then cuts $ H_{g}^{+} $ in a
bi-angle with inner angles $ \alpha \in \;]0,\pi[ $ at $ p $
and $ \beta \in \;]0,\pi[ $ at $ q $. Let $ \Gamma_{g}^{+} $
denote the oriented part of $ \Gamma $ in $ H_{g}^{+} $. The
Gauss/Bonnet theorem for this bi-angle says
\BE
\alpha+\beta = \int_{\Gamma_{g}^{+}} \kappa \; ds,
\label{eq_3.6}
\EE
hence
\BE
\alpha+\beta \leq \text{max}_{\,\Gamma_{g}^{+}}\, \kappa\cdot 
L(\Gamma_{g}^{+}).
\label{eq_3.7}
\EE
In particular, the bi-angle is acute (i.e. $ \alpha, \beta $
both less than $ \pi/2 $) if $ L(\Gamma_{g}^{+}) $ is small
enough.

Observe that, for an acute bi-angle, the tangents of 
$ \Gamma $ at $ p, q $ intersect in a point $ r $ of 
$ H_{g}^{+} $ and that $ \Gamma_{g}^{+} $ is contained in
the triangle with vertices $ p,q,r $. Furthermore, the arc 
$ \Gamma_{g}^{+} $ lies `schlicht' over the side $ [p,q] $, 
i.e. the straight lines orthogonally cutting the points 
$ z $ of $ [p,q] $ hit $ \Gamma_{g}^{+} $ exactly once, say
at $ h(z) $. The distance $ d(z,h(z)) $ is less than the
height $ h_{0} $ of that triangle over $ [p,q] $. By
elementary trigonometry, 
$ h_{0} = d(p,q)/(\cot\,\alpha+\cot\,\beta) $.

\bigskip
\BT[on approximation]
\label{lem_3.4}
~\\
Let $ \Gamma $ be an oriented closed and strictly convex 
$ C^{2} $-curve in $ \R^{2} $, and let $ k $ be an upper
bound of its curvature function $ \kappa: \Gamma \to \R $.
Let $ P $ be a polygon inscribed in $ \Gamma $ and 
$ \lambda $ an upper bound of the arclength on $ \Gamma $
between all consecutive vertices of $ P $ such that
\BE
\lambda < \frac{\pi}{2k}.
\label{eq_3.8}
\EE
Then we have for the minimax-invariants of $ \Gamma $ and 
$ P $
\BE
\delta(P) \leq \delta(\Gamma) \leq 
\delta(P)+\lambda\cdot\tan(k\lambda).
\label{eq_3.9}
\EE
Thus, if $ \lambda $ tends to $ 0 $, the minimax-invariant 
$ \delta(P) $ tends to $ \delta(\Gamma) $ with an error term 
of order $ \mathcal{O}(\lambda^{2}) $. 
\ET

\smallskip
\TI{Proof.} 
The vertices of the polygon $ P $ on $ \Gamma $ shall follow
each other according to the orientation of $ \Gamma $ ($ P $
itself is convex with interior in $ \Omega $). By
\eqref{eq_3.7}, \eqref{eq_3.8}, all bi-angles cut out by
the edges of $ P $ on $ C $ are acute having inner angles 
$ \leq k \lambda $. Let $ \psi: \Gamma \to P $ be the metric
projection onto $ P $, restricted to $ \Gamma $, and 
$ \vhi: P \to \Gamma $ the map whose restriction to any edge
$ [p,q] $ of $ P $ is the map $ h $, defined before. In
fact, $ \vhi $ and $ \psi $ are inverse to each other. It is
well known that $ \psi $ is Lipschitz-continuous with
Lipschitz-constant $ 1 $ (e.g. \TI{Phelps} [1957], Thm. 5.1),
hence $ d(p,q) \leq d(\vhi(p),\vhi(q)) $ for all 
$ p,q \in \Gamma $. On the other hand 
$ d(p,\vhi(p)) \leq \lambda /(2 \cot\,k \lambda) $. Thus,
Lemma \ref{lem_3.2} applies and yields \eqref{eq_3.9}.
\QED


\section{\hspace{-1em}. 
The Invariant \texorpdfstring{$ \pmb{\delta} $}{deltaX} for Polygons}
\label{dpoly_4}
\markright{\ref{dpoly_4}. The Invariant $ \delta $ for Polygons}


By the approximation theorem \ref{lem_3.4} it suffices to prove the desired inequality \eqref{eq_7} for convex polygons. In this case, the invariant $ \delta $ can be calculated in finitely many steps. The following observations will show this. Denote by 
$ P = \Gamma $ a convex polygon in $ \R^{2} $, by $ \Omega \neq \vac $ its inner region, and by $ C := \Omega \cup \Gamma $ its closure. For $ u, v \in \R^2 $ we denote by $ [u,v] $ the closed segment from $ u $ to $ v $, inverting the brackets at an endpoint which does not belong to the segment, e.g. 
$ \left[u,v\right[ := [u,v] \setminus \{v\} $. 

\BL
\label{lem_4_1}
Given a point $ x \in  \R^{2} $, the following holds true: If a point  $ y \in P $ has maximal distance from $ x $, then $ y $ is necessarily a vertex of $ P $.
\EL

\TI{Proof.} 
Of course, $ y \neq x $. If $ y $ would lie in the interior of an edge $ E $ then all points of $ P $ would lie in the closed disk 
$ \overline{\B}(x,r) $ with radius $ r:=\abs{x-y} $. On the other hand, the edge $ E $ is tangent to the circle $ S(x,r) $, by the maximality condition. This contains a contradiction, because a circle tangent minus the point of contact is always outside the circle.
\QED 

Henceforth, a point $ y \in P $ in maximal distance to an 
$ x \in \R^2 $ will be called a \TI{farthest} point from $ x $.
For given $ x \in  \R^{2} $ there may exist more than one 
farthest vertex from $ P $. This happens iff $ x $ lies on a perpendicular bisector of two vertices. The set of all perpendicular bisectors of any two vertices of $ P $ splits  
 $ \R^{2} $ in certain domains. More precisely: The complement of all these perpendicular bisectors is an open and dense subset of 
$ \R^{2} $ with finitely many connected components. These components arise as the possible intersections of open half spaces of all perpendicular bisectors (as far those intersections are nonvoid). For $ N $ perpendicular bisectors we thus have at most 
$ 2^{N} $ such connected components. 

If $ x \in  \R^{2} $ owes exactly one vertex $ y \in P $ as farthest point from $ x $ then all vertices $ \neq y $ have a smaller distance to $ x $. By the continuity of the distance function, this is retained when $ x $ varies in a small neighborhood. So the set of the $ x \in  \R^{2} $ which have $ y $ as a unique point of maximal distance on $ P $ is open. This implies that all points in one of the connected components from above have the same vertex as a point of maximal distance. 

If $ x $ varies only in $ P $ then this holds accordingly: The intersection points of all perpendicular bisectors with $ P $ form a partition of the polygon $ P $, while $ P $ is homeomorphic to 
the unit circle $ \SS^{1} $. This means: If $ M_{1},\ldots,M_{K} $ are these intersection points on $ P $, sorted in the positive mathematical sense along $ P $, then the complement 
$ P \setminus \{M_{1},\ldots,M_{K}\} $ is composed of the open and connected sets 
$ S_{1} := M_{1}M_{2} $,\ldots,$ S_{K-1} := M_{K-1}M_{K} $, 
$ S_{K} := M_{K}M_{1} $,
and the points of each such set all have the same uniquely determined farthest vertex. Here, for 
$ M \neq M' $ in $ P $, denotes $ MM' $ the part of $ P $ which traces from $ M $ to $ M' $ in the chosen orientation of $ P $ (without $ M $ and $ M' $). In the sequel, the points 
$ M_{1},\ldots,M_{K} $ will be called the \TI{section points}, and the sets $ S_{1} $,\ldots,$ S_{K} $ will be called the 
\TI{sections} of $ P $. For a section $ S $ we denote by  $ y(S) $ the unique farthest vertex, common to all $ x \in S $. 

It may very well happen that a section $ S $ contains vertices of 
$ P $. So, a section is not necessarily part of an edge. If one refines $ P $ by adding the vertices, then there arise new sections
which then lie on edges of $ P $. These new sections shall be called \TI{refined sections}. The points of a refined section 
$ S' $ then still have the same unique farthest vertex, namely the one, which belongs to the section $ S $ in which the subsection $ S' $ is contained. We then write $ y(S') := y(S) $.

The above function $ \mu $ may be determined on any section 
$ S $: One has  $ \mu(x) = d(x,y(S)) $ for all $ x \in S $ since 
$ d(x,q) \leq d(x,y(S)) $ for all $ q \in P $. By the continuity of 
$ \mu $ and $ d $ we still have
\BE
\label{eq_100108_7_*}
\mu(x) = d(x,y(S)) \quad \text{for $ x \in \overline{S} $}.
\EE
Herein, $ \overline{S} $ denotes the closure of $ S $, i.e. 
$ S $ together with its two endpoints. This implies that $ y(S) $ is still a farthest point from the two endpoints. However, such an endpoint owes at least one more farthest vertex 
$ \neq y(S) $ because it lies on a perpendicular bisector. 

The minimum of the function $ \mu $ on a (refined) section $ S $ equals the minimum of the function $ x \mapsto d(x,y(S)) $ for 
$ x \in S $ where $ y(S) $ does not depend on these $ x $. If, generally, the distance of a point $ z $ from a set $ T $ is denoted by
\BE
\label{dist_4}
d_{T}(z) := \Min \{d(p,z) \mid p \in T\}
\EE
then consequently
\BE
\label{minmu_4}
\min\nolimits_{S}\; \mu = \min\nolimits_{\overline{S}}\;\mu 
= d_{S}(y(S)) = d_{\overline{S}}(y(S)).
\EE
In case of a refined section $ S $ the set $ \overline{S} $ is a closed segment in $ \R^{2} $ whose distance from $ y(S) $ can be determined elementarily. Namely, if $ a $, $ b $ are the endpoints of the open segment $ S $, and if $ z \in \R^{2} $ is arbitrarily given then the nearest point $ z^{\Ast} $ for $ z $ on 
 $ \overline{S} $ is given by
 \BE
\label{eq_030908_4_1}
z^{\Ast} =
\begin{cases}
a+t(b-a) & \text{for $ 0 < t < 1 $} \\
a        & \text{for $ t \leq 0 $} \\
b        & \text{for $ t \geq 1 $}
\end{cases}
\qquad
\text{with}
\qquad
t := \frac{\Dot{z-a}{b-a}}{\abs{b-a}^{2}}.
\EE
The cases depend from whether the orthogonal projection of $ z $ onto the straight line spanned by $ S $ belongs to 
$ S $ (first case) or not (second and third case). In the first case, $ z^{\Ast} $ is the footpoint and the formula follows from the requirement that $ z-z^{\Ast} $ is orthogonal to $ b-a $. In the second and third case, $ z^{\Ast} $ is the nearest point 
for $ z $ on $ S $. From this, the distance of $ z $ to $ S $ is 
calculated as $ d_{S}(z) = \abs{z-z^{\Ast}} $, explicitly
\BE
\label{eq_100108_10_2}
d_{S}(z) = 
\begin{cases}
\ds\frac{\abs{(z-a\wedge(b-a)}}{\abs{b-a}} 
& \text{for $ 0 < t < 1 $} \\
\abs{z-a} & \text{for $ t \leq 0 $} \\
\abs{z-b} & \text{for $ t \geq 1 $}.
\end{cases}
\EE
where
\BE
\label{eq_100108_10_2a}
\abs{u \wedge v} := 
\sqrt{\abs{u}^{2} \cdot \abs{v}^{2}-\Dot{u}{v}^{2}}.
\EE
\TI{Remark.}
In \eqref{eq_030908_4_1}, the formula for  $ z^{\Ast} $ in the first case can also expressed by
\BE
\label{eq_030908_5}
z^{\Ast} = 
\frac{\Dot{z-b}{a-b}a-\Dot{z-a}{a-b}b}{\abs{a-b}^{2}}.
\EE
The global minimum $ \delta $ of the function $ \mu $ is now computable from these finitely many dates:

\BL
\label{lem__100108_10}
If $ S_{1} $,\ldots,$ S_{K} $ are the sections of $ P $ then
\BE
\label{eq_100108_10_1}
\delta(P) = 
\Min 
\{
d_{\overline{S}_{1}}(y(S_{1})),\ldots,
d_{\overline{S}_{K}}(y(S_{K})).
\}
\EE
The same formula holds if the $ S_{k} $ run through all 
\underline{refined} sections of $ P $ (then with different $ K $).
\EL

\TI{Proof.}
Generally one can say:

\TI{%
Let $ X = X_{1} \cup \cdots \cup X_{K} $ be a finite covering of a set $ X $ und $ f: X \to \R $ be a real function. If there exist the minima of all restrictions $ f|X_{k} $ in 
$ \R $ then there also exists the minimum of $ f $ in $ \R $,
and one has}
\BE
\label{eq_100108_8_1}
\Min f = 
\Min
\{
\Min f|X_{1},\ldots,\Min f|X_{K}
\}.
\EE
This implies the assertion.
\QED

For a convex polygon, any distinguished chord ends in a vertex 
(Lemma \ref{lem_4_1}).

Eqn.  \eqref{eq_100108_10_1} allows, in interaction with Eqn. 
\eqref{eq_100108_10_2}, the calculation of the minimax invariant 
$ \delta $ in a finite number of steps. If one also calculates, in each case, the nearest points $ y(S_{k}) $ on the sections by 
\eqref{eq_030908_4_1} or \eqref{eq_030908_5} then all distinguished chords can accordingly be obtained in finitely many steps.


\section{\hspace{-1em}. 
The Algorithm for Polygons}
\label{alg_4a}
\markright{\ref{alg_4a}. The Algorithm for Polygons}


Here we shall justify and describe in detail the above-mentioned algorithm. Given is only the list of vertices of the polygon $ P $ in the order along $ P $, counted mathematically positive. The aim is the determination of the invariant $ \delta(P) $ and of the set of all distinguished chords. In the notation, the dependence on 
$ P $ will be suppressed, so instead of $ \delta(P) $ we write 
$ \delta $, etc.

The sorting of all points of interest on $ P $ follows the   \TI{arclength}, measured from the first entry of the list in the counter clockwise manner. In particular, all the arclength parameters of the vertices can be calculated. The arclength parameter of the last entry (identical with the first) is then the perimeter $ L $ of $ P $.

\medskip
\TB{\TI{Calculation of $ \delta $}}

This takes place with Eqn. \eqref{eq_100108_10_1}.

First of all, the perpendicular bisectors of \TI{all} pairs of vertices must be cut with $ P $. (The perpendicular bisectors of the edges don't suffice.) Initially, the intersection points don't appear in the correct order along $ P $. However, the arclength parameters of these points can be determined edge-wise.
For a better performance one may consider that one of the perpendicular bisectors cuts the polygon in exactly two points. If the perpendicular bisector belongs to an edge then one of these two points the midpoint of the edge. In addition, one may arrange the section points in the order of increasing arclength. With this, the sections can be treated successively, by adding accordingly the arclength parameters of the vertices. In fact, one must work with the refined sections, since the formulae 
\eqref{eq_030908_4_1} to \eqref{eq_030908_5} refer to a segment but not to an open polygon chain. (It may happen that a section contains points of multiple edges.) 

For each refined section, then the farthest vertex and its shortest distance to the section can be calculated. Subsequently, 
$ \delta $ can be calculated by a finite minimum according to 
\eqref{eq_100108_10_1}.

\medskip
\TB{\TI{Calculation of the distinguished chords}}

Let $ S $ be a refined section which realizes the minimum 
\eqref{eq_100108_10_1}, so with $ \delta = d_{S}(y(S)) $. The two following lemmata describe how to obtain from such a refined section the distinguished chords.

\BL
\label{lem_4a_1}
Let $ S $ be a refined section with $ \delta = d_{S}(y(S)) $.

If  $ z^{\Ast} \in \overline{S} $ is the nearest point for 
$ y(S) $, hence $ \delta = d_{S}(y(S)) = d(z^{\Ast},y(S)) $ then the segment $ [z^{\Ast},y(S)] $ is a distinguished chord of $ P $.
\EL

\TI{Proof.}
One has to show
\BE
\label{eq_4a_1_1}
d(z^{\Ast},y(S)) \geq d(z^{\Ast},q) 
\quad \Forall q \in P.
\EE
For this it suffices to show
\BE
\label{eq_4a_1_2}
d(z^{\Ast},y(S)) = \max\nolimits_{q\in P} d(z^{\Ast},q).
\EE
Setting in \eqref{eq_100108_7_*}: $ x = z^{\Ast} $ implies
$ \mu(z^{\Ast}) = d(z^{\Ast},y(S)) $. By 
$ \mu(z^{\Ast}) = \max\nolimits_{q\in P} d(z^{\Ast},q) $, this expresses the assertion.
\QED

The farthest vertex $ y(S) $ from the \TI{inner} points of an arbitrary (refined) section $ S $ is always uniquely determined.
This does not hold necessarily for the two ends of $ S $. Let, for instance, $ A $ be the starting point of $ S $. Admittedly, 
$ y(S) $ is also the farthest vertex from $ A $, but $ A $ may owe further farthest vertices.

\BL
\label{lem_4a_2}
Let $ S $ be a refined section with $ \delta = d_{S}(y(S)) $.

If, for an initial point $ a $ of $ S $ and a vertex 
$ \Tilde{a} \neq a $ one has
\BE
\label{eq_4a_2_1}
\delta = d(a,y(S)) = d(a,\Tilde{a})
\EE
then $ [a,\Tilde{a}] $ is also a distinguished chord of $ P $.
\EL

\TI{Proof.}
The segment $ [a,y(S)] $ is a distinguished chord, by Lemma \ref{lem_4a_1}, so
\BE
\label{eq_4a_2_2}
\delta = d(a,y(S)) \geq d(a,q) \quad \Forall q \in P.
\EE
The same inequality holds, by \eqref{eq_4a_2_1}, if  $ d(a,y(S)) $ is replaced by $ d(a,\Tilde{a}) $. Thus, $ [a,\Tilde{a}] $ is a distinguished chord.
\QED

Under the hypothesizes just made, $ a $ necessarily is a section point because $ y(S) $ and $ \Tilde{a} $ have the same distance from $ a $, so  $ a $ lies on the perpendicular bisector of the chord $ [y(S),\Tilde{a}] $. 

Obviously, the algorithm can be shaped in such a way that all distinguished chords which occur by the Lemmata \ref{lem_4a_1} und \ref{lem_4a_2} are generated in finitely many steps. In fact, 
\TI{all} existing distinguished bisectors can be obtained in this way:

\BL
\label{lem_4a_3}
Any distinguishes chord $ [p_{0},q_{0}] $ of the polygon $ P $ 
arises by the assumptions of Lemma \ref{lem_4a_1} or Lemma \ref{lem_4a_2}.
\EL

\TI{Proof.}
At any rate, there exists a refined section $ S $ such that 
$ p_{0} $ lies in $ S $ or is the initial point of $ S $. By assumption one has
\BE
\label{eq_4a_3_1}
\delta = d(p_{0},q_{0}) \geq d(p_{0},q) 
\quad \Forall q \in P.
\EE
In particular, $ q_{0} $ is a farthest vertex from $ p_{0} $. 
The following cases may occur:
\begin{DES1}{\TB{(ii)}}
\item[(i)]
\TI{$ q_{0} $ is also a farthest vertex for the points of 
$ S $:}
Then $ q_{0} = y(S) $, because such a farthest vertex is uniquely determined. First, one has to show: 
\BE
\label{eq_4a_3_2}
d_{T}(y(T)) \geq d_{S}(y(S)) \quad 
\text{for all refined segments $ T $}:
\EE
This follows from 
\BE
\label{eq_4a_3_3}
d_{T}(y(T)) \geq \delta = d(p_{0},q_{0}) = d(p_{0},y(S)) 
\geq d_{\overline{S}}(y(S)) = d_{S}(y(S)).
\EE
Specifically for $ T = S $, here appears everywhere the equality sign: $ d(p_{0},q_{0}) = d_{S}(y(S)) $. So $ p_{0} $ is the nearest point for  $ y(S) $ in $ \overline{S} $. Hence for 
$ z^{\Ast} := p_{0} $, all assumptions of Lemma \ref{lem_4a_1} are 
fulfilled.
\item[(ii)]
\TI{$ q_{0} $ is \underline{not} a farthest vertex for the points of $ S $:}
Then  $ p_{0} \notin S $, yet $  p_{0} =: a $ is the initial point of $ S $. By the assumption \eqref{eq_4a_3_1} one has:
$ \delta(a,q_{0}) \geq d(a,q) $ for all  $ q \in P $. Thus, 
$ q_{0} $ is a farthest point from $ a $, but $ q_{0} \neq y(S) $.

One has $ \delta = d(a,q_{0}) = d_{S}(y(S)) $ because $ q_{0} $ as well as $ y(S) $ are farthest vertices from $ a $. So, for 
$ S, a $ and $ \Tilde{a} := q_{0} $, the assumptions of Lemma 
\ref{lem_4a_2} are satisfied. 
\QED
\end{DES1}

Lemma \ref{lem_4a_3} justifies the algorithm which rests on the Lemmata \ref{lem_4a_1} and \ref{lem_4a_2}.

In the following pictures, the endpoints of the sections will be marked by little crosses and the distinguished chords as dashed segments.


\section{\hspace{-1em}. 
Short polygons (\texorpdfstring{$ n \in \{3,4\} $}{n34})}
\label{udr_5}
\markright{\ref{udr_5}. Short polygons ($ n \in \{3,4\} $)}


The general polygon case with $ n $ edges is completely open. Though many experiments seem to support the `isoperimetric inequality' \eqref{eq_7}, we can prove it only for triangles and offer a concrete conjecture for quadruples, the latter with a surprising figure as a `magic quadrangle' in the equality case.

Restricting oneself to the pure polygon case one may expect that
the `isoperimetric quotient' $ L/\delta $ assumes a value bigger than $ \pi $ but with limit $ \pi $ for $ n \to \infty $. 

\medskip
\TB{\TI{The case of triangles}}

\BT
\label{thm_undr}
For each triangle $ \Delta $ in $ \R^{2} $ one has
\BE
\label{eq_230808_65_1}
L(\Delta) \geq 2\sqrt{3}\cdot \delta(\Delta),
\EE
with equality exactly for equilateral triangles.
\ET

The quotient $ L/\delta $ is invariant under similarities. On may use this and other restrictions in order to build a 
`moduli space' for all convex polygons of a fixed edge number 
$ n $, thus reducing the problem to a study of $ L/\delta $ on this moduli space. For triangles we shall follow a principle which is also sensible in higher cases, namely to single out a diameter and normalize it to a fixed size.

\medskip
\TB{\TI{Diameter in the general polygon case}}

\BL
\label{lem_diam}
A diameter $ D $ of a convex polygon $ P $ always connects two 
vertices. 

Denote by $ H_{1} $ and $ H_{2} $ the two open half planes of the straight line spanned by $ D $ and by $ D^{\Circ} $ the relative interior of $ D $.

If $ P \cap H_{i} \neq \vac $ for an $ i \in \{1,2\} $ then 
$ P \cap \overline{H}_{i} \setminus D^{\Circ} $ lies schlicht over 
$ D $ in the direction orthogonal to $ D $. If  
$ P \cap H_{i} = \vac $ then $ P \cap \overline{H}_{i} $ is an edge of $ P $ and equal to $ D $.
\EL 

\TI{Proof.}
The first part follows directly from Lemma \ref{lem_4_1} because each endpoint of a diameter is the farthest point from the other endpoint.

($ {\Ast} $)
Without restriction, let the endpoints of the diameter $ D $ be  
$ A := (-1,0) $ and $ B := (1,0) $ and assume that  $ H_{1} $ is the upper half plane and $ H_{2} $ the lower half plane of 
$ \R^{2} $.

Each of the four edges through $ A $ or $ B $ makes an angle 
$ < \pi/2 $ with $ D $ because, otherwise, there resulted a contradiction to the maximal length of $ D $. If one of these edges forms an angle $ 0 $ with $ D $ then $ D $ is an edge of 
$ P $ and $ P $ lies completely in $ \overline{H}_{1} $ or in 
 $ \overline{H}_{2} $.

If e.g. $ P \cap H_{1} \neq \vac $ then always one of the four edges through $ A $, $ B $ is directed into $ H_{1} $ (while the two others don't this). The part of $ P $ in $ \overline{H}_{1} $ forms, together with $ D $ a convex polygon 
$ P_{1} := (P \cap \overline{H}_{1}) \cup D $ with $ D $ as an edge and diameter. Then the two other edges of $ P_{1} $ through 
$ A $, $ B $ make with $ D $ an angle between $ 0 $ and $ \pi/2 $.
In this situation, a straight line orthogonal to $ D $, 
$ te_{1}+\R e_{2} $, doesn't cut the polygon $ P_{1} $ if 
$ t < -1 $ oder $ t > 1 $, and except in $ te_{1} $ exactly once in a point $ te_{1}+f(t)e_{2} $ with $ f(t) \geq 0 $ if 
$ -1 \leq t \leq 1 $. Otherwise, in case $ \abs{t} > 1 $, there were a contradiction to the diameter property of $ D $. In case 
$ \abs{t} < 1 $, the aforementioned straight line cannot have a whole segment in common with $ P_{1} $ since, then, it were a support line of $ P_{1} $ such that either $ A $ or $ B $ would not belong to $ P_{1} $. In case $ t = \mp 1 $, obviously $ A $, resp. $ B $ is the sole common point of the straight line with 
$ P_{1} $. So the convex region $ C_{1} $ bounded by $ P_{1} $ is given by 
\BE
\label{eq_lem_diam_1}
C_{1} = 
\{(t,s) \mid -1 \leq t \leq 1, \; 0 \leq s \leq f(t) \}
\EE
with a function $ f: [-1,1] \to \R $ with $ f(-1) = f(1) = 0 $ and $ f(t) \geq 0 $ for $ \abs{t} < 1 $. In fact, $ f $ is a concave function since $ C_{1} $ is convex and $ C_{1} $ is the subgraph of $ f $. By $ f(-1) = f(1) = 0 $ this implies $ f(t) > 0 $ for 
$ \abs{t} < 1 $.

Projecting the edges of $ P_{1} $ onto $ D $ one sees that the function $ f $ is piecewise affine. 

If, e.g., $ P \cap H_{2} = \vac $ then $ P \cap H_{1} \neq \vac $ and the same arguments as above for $ P_{1} $ apply for $ P $ itself.

We still record:

\BR
\label{rem_added}
In case $ P \cap H_{1} \neq \vac $ there holds, under 
the normalization ($ \Ast $): 
\BE
\label{lem_diamz}
P \cap \overline{H}_{1} = 
\{(t,f(t)) \mid -1 \leq t \leq 1 \} \cup D,
\EE
with a concave and piecewise affine function
$ f: [-1,1] \to \R $, where $ f(-1) = f(1) = 0 $ and $ f(t) > 0 $ for $ \abs{t} < 1 $. 
\QED
\ER

\TB{\TI{Proof of Theorem \ref{thm_undr}}}

In the triangle case things shall be arranged in the above manner.
A diameter is necessarily one of the longest edges, so 
the segment from $ A = (-1,0) $ to $ B = (1,0) $ is an edge of maximal length. The edge of the middle size may start in $ A $, the shortest edge in $ B $. The third vertex $ C = (x,y) $ then lies in the quadrant $ \{(x,y) \mid x \geq 0,\; y \geq 0\} $. The moduli space sounds here
\BE
\label{eq_modr_1}
\mathcal{M} := 
\{(x,y) \mid x \geq 0,\;y \geq 0,\; (x+1)^{2}+y^{2}\leq 4\},
\EE
where, admittedly, the segment from $ (0,0) $ to 
$ B $ ($ \!\! \iff y = 0 $) corresponds to degenerate triangles
($ \Omega = \vac $). In the sequel, it is assumed that $ y > 0 $.

The centers of the edges $ AB $, $ BC $, $ CA $ are
\BE
\label{eq_modr_2}
M := (0,0), \qquad 
M_{+} := \frac{1}{2}(x+1,y), \qquad
M_{-} := \frac{1}{2}(x-1,y).
\EE
The section points of the perpendicular bisectors through 
$ M $, $ M_{+}$, $ M_{-} $ with $ P $ are 
\BE
\label{eq_modr_3}
S = (0,h), \qquad
S_{+} = (x_{+},0), \qquad
S_{-} = (x_{-},0),
\EE
with
\BE
\label{eq_modr_4}
h := \frac{y}{1+x}, \qquad
x_{+} := \frac{1}{2}\;\frac{1-(x^2+y^2)}{1-x}, \qquad
x_{-} := \frac{1}{2}\;\frac{(x^2+y^2)-1}{1+x}.
\EE
From this one reads that $ M_{+} $ always lies right of the 
$ y $-axis and $ M_{-} $ left of it and, further, that $ S_{+} $ and $ S_{-} $ are situated on the same resp. different sides of the $ y $-axis according to $ x^2+y^2 < 1 $ resp. $ x^2+y^2 > 1 $.
Also, the angle at $ C $ follows these inequalities: In the first case it is obtuse, in the second it is sharp. This results from the calculation of the scalar product 
$ \Dot{A-C}{B-C} = x^2+y^2-1 $. Cf. the following pictures:

\vspace{-2ex}
\begin{center}
\begin{minipage}{15.8cm}
\bild{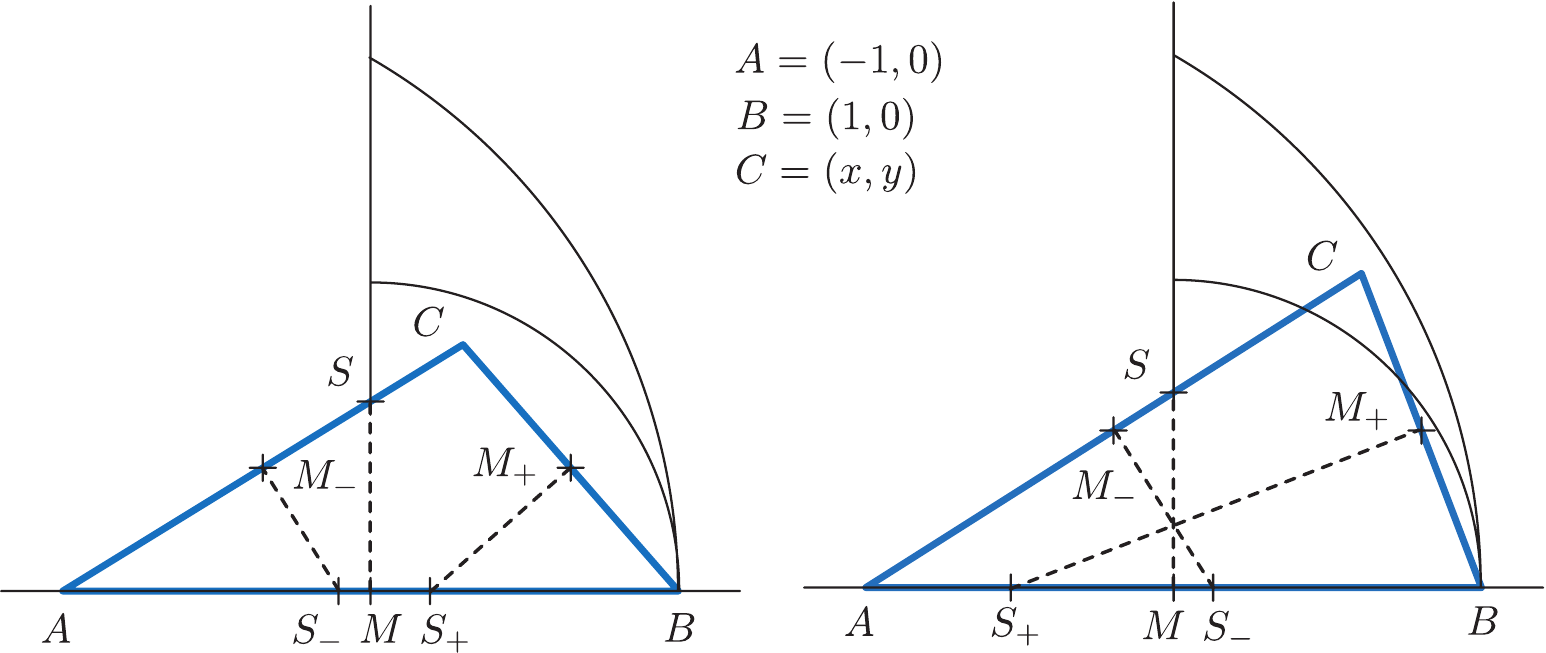}
\end{minipage}
\end{center}

In case $ x^2+y^2 = 1 $ (right angle at $ C $) nothing has to be changed w.r.t. $ M_{+} $ and $ M_{-} $, but 
$ S_{+} $ und $ S_{-} $ coincide with $ M $.

In the sequel we will refer to this situation as far as it is geometrically evident.

One realizes that the unit circle decomposes the moduli space 
$ \mathcal{M} $ into two parts which behave differently. In fact, the differences will become even more drastic with regard to the invariant $ \delta $.

\medskip
\TB{\TI{Subcase $ x^2+y^2 \leq 1 $}}

Under this assumption the behaviour is pretty simple and also uniform. 

\BL
\label{lem_untb}
If, for $ C \in \mathcal{M} $, one has $ x^2+y^2 \leq 1 $ then 
$ \delta = 1 $, hence by $ L > 4 $:
$$
\frac{L}{\delta} > 4.
$$
\EL

For those $ C \in \mathcal{M} $ the inequality \eqref{eq_230808_65_1} is strict, so nothing has to be added.

\TI{Proof of \ref{lem_untb}.} 
First of all, the case $ x^2+y^2 < 1 $ will be treated explicitly.
The decomposition of $ P $ into sections looks as follows:
 
\vspace{-2ex}
\begin{center}
\begin{minipage}{10.6cm}
\bild{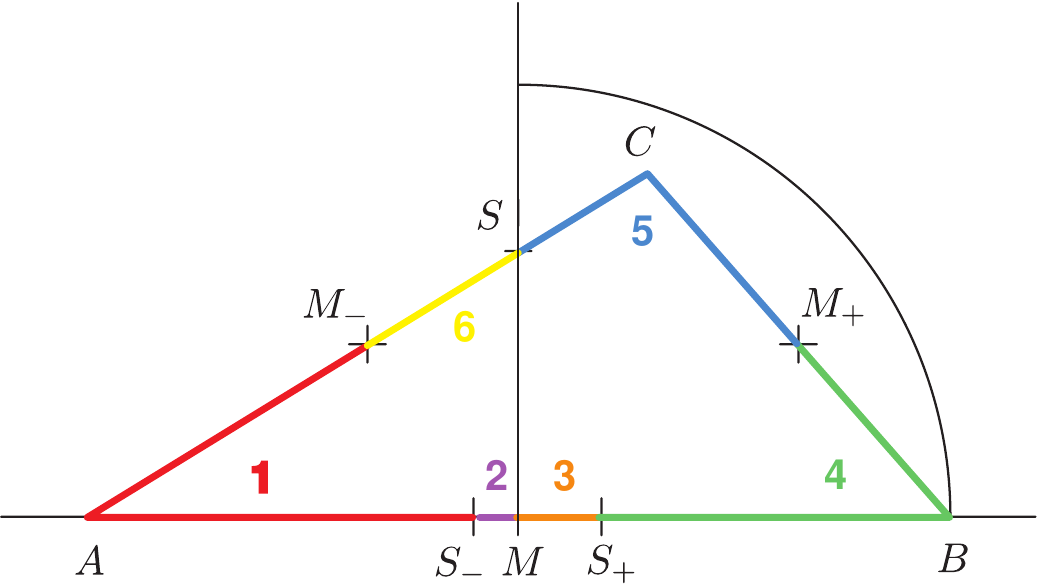}
\end{minipage}
\end{center}

In particular, some of the sections `round the corner' (namely those with numbers $ 1 $, $ 4 $, and $ 5 $). This means that they contain vertices which are not section points themselves. If the vertices would be added to the decomposition then from the six sections there would arise nine refined sections. However, these nine sections have no impact on the following discussion.

For short, write $ y_{k} := y(S_{k}) $ and 
$ d_{k} := d_{S_{k}}(y(S_{k})) $, and denote by $ z_{k}^{\Ast} $ a nearest point to $ y_{k} $ on $ S_{k} $. These objects are determined according to the following table:
\BE
\label{eq_untb_1}
\begin{array}{l|c|c|c}
S_{k} & y_{k} & z_{k}^{\Ast}              & d_{k} \\ \hline
S_{1} & B     & M_{-} \text{ or } S_{-} & > 1 \\[0.5ex]
S_{2} & B     & M                         & 1   \\[0.5ex]
S_{3} & A     & M                         & 1   \\[0.5ex]
S_{4} & A     & M_{+} \text{ or } S_{+} & > 1 \\[0.5ex]
S_{5} & A     & M_{+} \text{ or } S     & > 1 \\[0.5ex]
S_{6} & B     & M_{-} \text{ or } S     & > 1
\end{array}
\EE
Most entries are immediately obvious. For example, for the section 
$ S_{1} $, one has $ y_{1} = B $ because $ A \in S_{1} $ and also 
$ z_{1}^{\Ast} $ is equal to $ M_{-} $ or to $ S_{-} $, so in any case $ d_{1} > 1 $. For the section $ S_{2} $ one has 
$ y_{2} = B $, as shown by comparison with $ d(M,C) = 1 $, and evidently $ z_{2}^{\Ast} = M $, hence $ d_{2} = 1 $. The following lines in the table are similarly simple. In some cases it is of some importance that the orthogonal footpoints of $ A $ und $ B $ on the opposite edges come to lie outside of these. \footnote{A more detailed determination of $ z_{k}^{\Ast} $ in the `or-cases' of the table \eqref{eq_untb_1} is possible but not really necessary here.}

By this, the minimal $ d_{k} $ is $ 1 $, thus $ \delta = 1 $. Distinguished chords by Lemma \ref{lem_4a_1} are those from 
$ M $ to $ A $ and to $ B $. To be added is another distinguished chord by Lemma \ref{lem_4a_2} from $ M $ to $ C $.

The case $ x^2+y^2 = 1 $ resolves itself in the same manner
or, also, by performing the limit $ x^2+y^2 \upto 1 $,  
applying Lemma \ref{lem_3.2}, part (i). 
Namely, a triangle with $ x^2+y^2 = 1 $ can be obtained in a
continuous manner from a triangle with $ x^2+y^2 < 1 $, e.g.
by a radial deformation with center $ 0 $.  Hereby, the
relation between the triangles can be kept bijectively and
uniformly small.

\medskip
\TB{\TI{Subcase $ x^2+y^2 > 1 $}}

Here, more possibilities are to be considered, and the discussion 
is clearly more complicated. The decomposition of $ P $ looks as follows:

\begin{center}
\begin{minipage}{10.4cm}
\bild{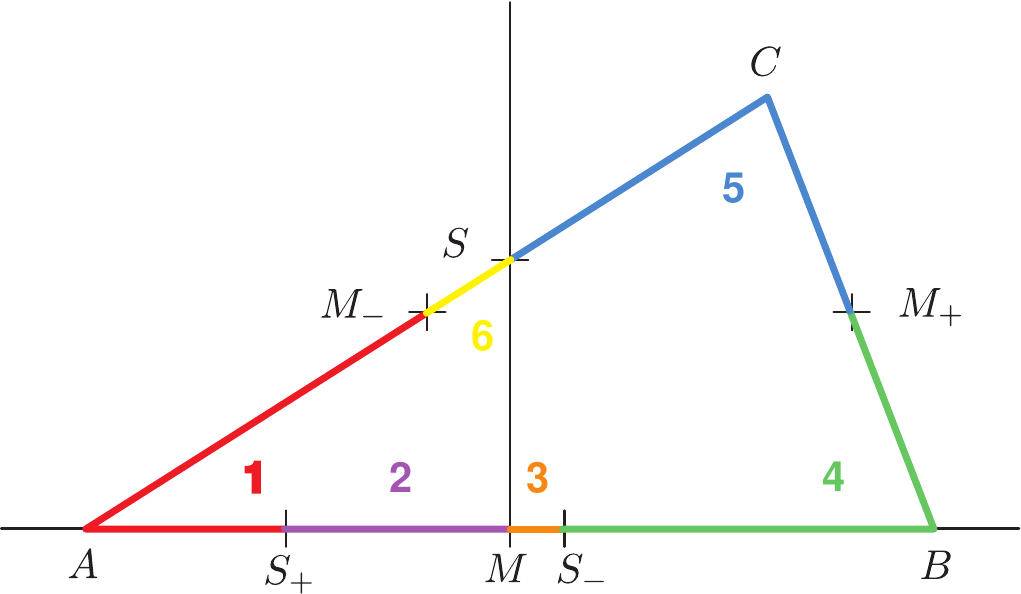}
\end{minipage}
\end{center}

In order that the sections look really that way (i.e. that the five points on  $ [A,B] $ and the four points on $ [A,C] $ are pairwise distinct) it is necessary to assume:
\BE
\label{eq_obb_1}
x > 0, \quad (x+1)^2+y^2 < 4.
\EE
One can confirm this constellation with the formulae \eqref{eq_modr_2} to \eqref{eq_modr_4}. By \eqref{eq_obb_1}, the two equilateral cases $ d(A,S) = d(B,S) $ resp. 
$ d(A,C) = d(A,B) = 2 $ are excluded (they will be discussed later on). 

The difference to the preceding case is the interchange in the position of $ S_{+} $ und $ S_{-} $ and the fact that a sharp angle now arises at $ C $.

It is suitable to calculate first the farthest points of the sections. Hereby, one can save some effort by regarding the following principle on the role of the perpendicular bisectors as discussed in section \ref{dpoly_4}:

\BL
\label{lem_obb_1}
If a point $ X \in P $ moves strictly monotonously towards a point 
which lies on exactly one perpendicular bisector (of two vertices 
$ P, Q $) and if the farthest vertex from $ X $ has been 
$ P $ before that it is equal to $ Q $ afterwards.
\QED
\EL

It is immediately clear that 
\BE
\label{eq_obb_2}
y_{1} = B, \quad y_{4} = A, \quad y_{5} = A,
\EE
since the sections with numbers $ 1,4,5 $ each contain a vertex. With the principle \ref{lem_obb_1} this implies 
\BE
\label{eq_obb_3}
y_{2} = C, \quad y_{3} = C, \quad y_{6} = B.
\EE
Now the distances $ d_{k} $, i.e. the footpoints $ z_{k}^{\Ast} $ of $ y_{k} $ on $ \overline{S}_{k} $ must be determined.  

\TI{For $ z_{1}^{\Ast} $:}

The footpoint of $ y_{1} = B $ on the straight line $ A \vee C $ lies outside the segment $ [A,M_{-}] $ because the angle 
$ \sphericalangle(AM_{-}B) $ is obtuse. This follows from a calculation which relies on Eqns. \eqref{eq_modr_2} to \eqref{eq_modr_4}:
\BE
\label{eq_obb_4}
\Dot{A-M_{-}}{B-M_{-}} = \frac{1}{4}((x-1)^2+y^2-4) < 0.
\EE
The last estimate is valid since otherwise one had $ x < 0 $, using \eqref{eq_obb_1}. 

Thus the foot of $ B $ on  $ \overline{S}_{1} $ equals $ M_{-} $ or $ S_{+} $. The decision for this results from the sign of the quantity
\BE
\label{eq_obb_5}
\vhi_{1}(x,y) := 
4\cdot(\abs{M_{-}-B}^2-\abs{S_{+}-B}^2) =
(1-x)^2((x-3)^2+y^2)-((x-1)^2+y^2)^2.
\EE
For $ x \in [0,1] $, the sign condition for $ \vhi_{1}(x,y) $ can be rearranged by solution for $ y $. This yields 
\BE
\label{eq_obb_6}
z_{1}^{\Ast} =
\begin{cases}
M_{-} & \quad \text{for $ y \geq \psi_{1}(x) $} \\
S_{+} & \quad \text{for $ y \leq \psi_{1}(x) $},
\end{cases}
\EE
where
\BE
\label{eq_obb_7}
\psi_{1}(x) := 
\sqrt{\frac{1-x}{2}} \sqrt{\sqrt{(9-x)^2-48}-(1-x)}.
\EE
In the equality case of \eqref{eq_obb_6}, both points $ M_{-} $ and $ S_{+} $ are footpoints of $ y_{1} $ on $ S_{1} $ (with same distances).

The graph of the function $ \psi_{1} $ over $ \left[0,1\right[ $ lies there between the circles $ \SS(M,1) $ and $ \SS(A,2) $ (with equality at $ x = 1 $); for one confirms that
\BE
\label{eq_obb_8}
\sqrt{1-x^2} < \psi_{1}(x) < \sqrt{4-(1+x)^2} \quad 
\text{for $ x \in \left[0,1\right[ $}.
\EE
\TI{For $ z_{2}^{\Ast} $:}

Obviously $ z_{2}^{\Ast} = M $ (consider the circle through 
$ M $ around $ y_{2} = C $).

\TI{To $ z_{3}^{\Ast} $:}

The footpoint of $ y_{3} = C $ on $ S_{3} $ is equal to  $ S_{-} $ or to the footpoint $ (x,0) $ of $ C $ on the straight line 
$ A \vee B $, according to the sign of  $ x-x_{-} $. One calculates  
\BE
\label{eq_obb_9}
x-x_{-} = \frac{(1+x)^{2}-y^{2}}{2(1+x)},
\EE
hence
\BE
\label{eq_obb_10}
z_{3}^{\Ast} =
\begin{cases}
S_{-} & \quad \text{for $ y \geq \psi_{3}(x) $} \\
(x,0) & \quad \text{for $ y \leq \psi_{3}(x) $},
\end{cases}
\qquad
\psi_{3}(x) := 1+x.
\EE

\TI{For $ z_{4}^{\Ast} $:}

As  $ z_{4}^{\Ast} $, only the points $ M_{+} $, $ S_{-} $ come in consideration. The foot of $ y_{4} = A $ on the segment 
$ [M_{+},B] $ is really equal to $ M_{+} $ because the angle $ \sphericalangle(AM_{+}B) $ is obtuse. This follows from the calculation
\BE
\label{eq_obb_11}
\Dot{A-M_{+}}{B-M_{+}} = 
\frac{1}{4}((1+x)^{2}+y^{2}-4) < 0.
\EE
In fact, the point $ S_{-} $ is always closer to $ S_{4} $ as 
$ M_{+} $. This rests on the inequalities
\BE
\label{eq_obb_12}
d(A,M_{+}) > 1+\frac{1}{2}(1+x) > x_{-}.
\EE
The first inequality expresses the fact that the orthogonal projection of the segment $ [A,M_{+}] $ onto the $ x $-axis shortens the segment. The second inequality reduces itself to 
$ y < \sqrt{4+4x} $ by conversion, and this is certainly true for 
$ y \leq \sqrt{4-(1+x)^{2}} $. Thus:
\BE
\label{eq_obb_13}
z_{4}^{\Ast} = S_{-}.
\EE

\TI{For $ z_{5}^{\Ast} $:}

Similar to the case $ z_{4}^{\Ast} $, the footpoint $ L_{5} $ of 
$ y_{5} = A $ on the straight line $ B \vee C $ now lies on the segment $ [M_{+},C] $. In the rectangular triangle $ A L_{5}C $, 
the leg $ [A, L_{5}] $ is shorter than the hypotenuse $ [A,C] $, so
\BE
\label{eq_obb_14}
z_{5}^{\Ast} = S.
\EE

\TI{For $ z_{6}^{\Ast} $:}

By  \eqref{eq_obb_4}, the angle $ \sphericalangle(BM_{-}S) $ is sharp, hence the footpoint of $ y_{6} = B $ on $ S_{6} $ is equal to 
$ S $ or to the footpoint $ L_{6} $ of $ B $ on the straight line $ C \vee A $. Decisive for this is the scalar product
\BE
\label{eq_obb_15}
\Dot{A-S}{B-S} = h^{2}-1 = \frac{(y-(1+x))(y+1+x)}{(1+x)^{2}}.
\EE
If the scalar product is $ \leq 0 $ then $ z_{6}^{\Ast} = S $; if it is $ \geq 0 $ then $ z_{6}^{\Ast} = L_{6} $ where by \eqref{eq_030908_5}:
\BE
\label{eq_obb_16}
L_{6} = 
\frac{\Dot{B-C}{A-C}A-\Dot{B-A}{A-C}C}{\abs{A-C}^{2}}.
\EE
The distance $ d(B,S) $ resp. $ d(B,L_{6}) $ itself is calculated as follows, using \eqref{eq_100108_10_2}
\BE
\label{eq_obb_17}
d_{6} :=
\begin{cases}
\ds \frac{2y}{\sqrt{(1+x)^{2}+y^{2}}} &
\quad \text{for $ y \geq \psi_{3}(x) $} \\[3ex]
\ds \frac{\sqrt{(1+x)^{2}+y^{2}}}{1+x} &
\quad \text{for $ y \leq \psi_{3}(x) $}
\end{cases}
\EE
with $ \psi_{3} $ as in \eqref{eq_obb_10}.

The region which is now of interest lies between the circle arcs 
$ \SS(M,1) $ and $ \SS(A,2) $ (always in the first quadrant) and is divided into four subregions I, II, III, IV, according to the cases \eqref{eq_obb_6} and \eqref{eq_obb_10}. Cf. the following picture

\begin{center}
\begin{minipage}{10.4cm}
\bild{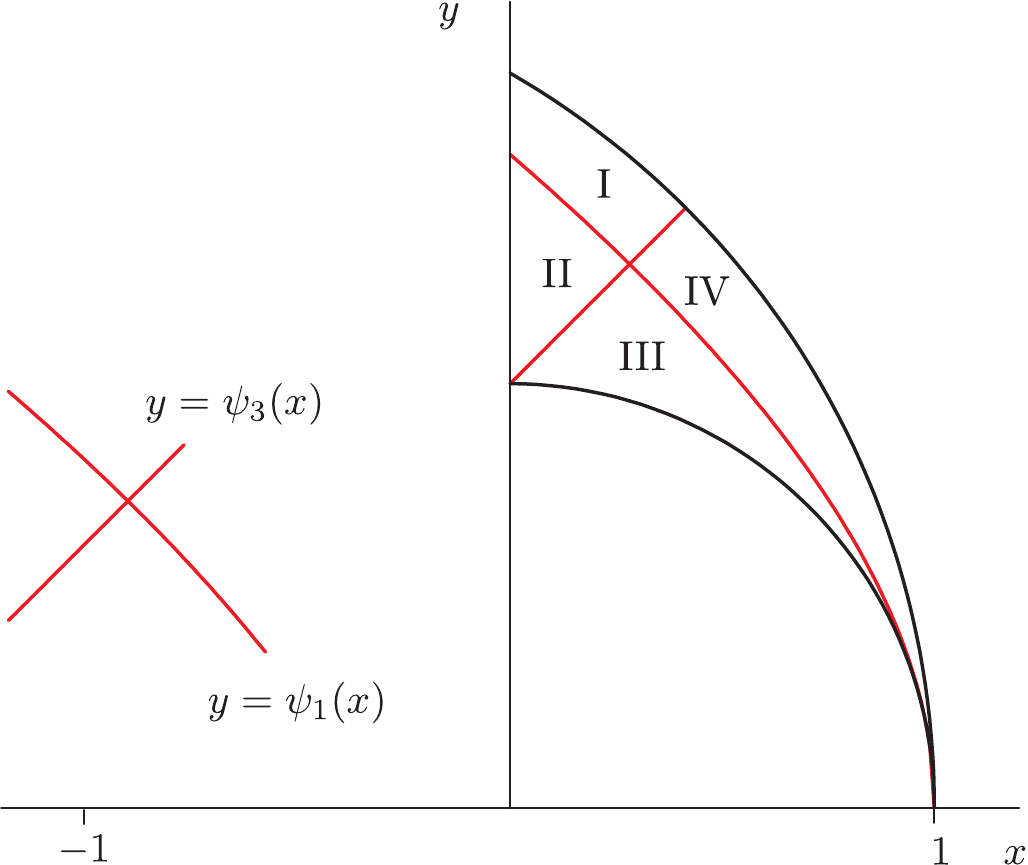}
\end{minipage}
\end{center}

The key at left indicates the equations describing the separating curves. The following features will justify this appearance.

The function $ \psi_{1} $ has the boundary values 
$ \psi_{1}(0) = \sqrt{\sqrt{33}-1}/\sqrt{2} = 1.5402\ldots $, 
$ \psi_{1}(1) = 0 $ and it decreases strictly monotonously in between. Namely, the derivative of $ \psi_{1} $ on 
$ \left]0,1\right[ $ equals, up to a positive factor, the following expression which can be easily estimated:
\BE
\label{eq_obb_17a}
(1-x)\cdot\sqrt{(x-9)^{2}-48}-((x-7)^{2}-28) < 
1\cdot\sqrt{16}-21 < 0.
\EE
The graph of the function $ \psi_{3} $ over $ [0,1] $ cuts the circle $ \SS(A,2) $ at the point  
$ (\sqrt{2}-1,\sqrt{2}) = (0.4142\ldots,1.4142\ldots) $. 
Since $ \psi_{1}(\sqrt{2}-1) = \sqrt{8\sqrt{2}-10} $, the two graphs of $ \psi_{1} $ and $ \psi_{3} $ over $ [0,\sqrt{2}-1] $ intersect themselves exactly once. The condition for this intersection leads to the cubic equation $ -3+12x-6x^2+4x^3+x^4 = 0 $ which has a unique real solution in $ [0,1] $ (computed numerically as) $ 0.2817\ldots $. So, the section point sounds 
$ (0.2817\ldots,1.2817\ldots) $.

Now, for each region, we must deduce the six distances 
$ d_{k} = d(y_{k},z_{k}^{\Ast}) $ from the above information on the foots $ z_{k}^{\Ast} $ and also their minima. The following table contains the results for the distances:  
\BE
\label{eq_obb_19}
\begin{array}{l|c|c|c|c}
d_{k} & \text{region I} & \text{region II} & \text{region 
III} & \text{region IV}\\ \hline 
d_{1} & \text{\rule[-2mm]{0mm}{10mm}}\ds\frac{1}{2}\sqrt{(3-x)^{2}+y^{2}} &
\ds\frac{1}{2}\;\frac{(1-x)^{2}+y^{2}}{1-x} &
\ds\frac{1}{2}\;\frac{(1-x)^{2}+y^{2}}{1-x} &
\ds\frac{1}{2}\sqrt{(3-x)^{2}+y^{2}} \\[2ex]
d_{2} & d_{2} &
d_{2} &
d_{2} &
d_{2} \\[2ex]
d_{3} & y &
y &
d_{4} &
d_{4} \\[2ex]
d_{4} & d_{4} &
d_{4} &
d_{4} &
d_{4} \\[1ex]
d_{6} & \ds\frac{2y}{\sqrt{(1+x)^{2}+y^{2}}} &
\ds\frac{2y}{\sqrt{(1+x)^{2}+y^{2}}} &
&
\end{array}
\EE
where
\BE
\label{eq_obb_20}
d_{2} = \sqrt{x^{2}+y^{2}}, \quad
d_{4} = \ds\frac{1}{2}\;\frac{(1+x)^{2}+y^{2}}{1+x}, \quad
d_{5} = \ds\frac{\sqrt{(1+x)^{2}+y^{2}}}{1+x}.
\EE
Obviously
$$ 
d_{5} = \sqrt{\frac{2}{1+x}}\sqrt{d_{4}} > d_{4}
$$
with equality on $ \SS(A,2) $). So, for the search of the minima, the $ d_{5} $ can be ignored. Therefore the line for $ d_{5} $ is omitted in the table, and similar for the two entries downright (which also contained $ d_{5} $). 

\TB{\TI{Intermezzo}}

In some cases, the handling of the inequalities is facilitated by passing from the cartesian coordinates $ x,y $ to the elliptic coordinates $ u,v $ defined by
\BE
\label{eq_ell_1}
\begin{aligned}
u &= 
\frac{1}{2}
\left(\sqrt{(1+x)^{2}+y^{2}}+
\sqrt{(1-x)^{2}+y^{2}}\right) \\
v &= 
\frac{1}{2}
\left(\sqrt{(1+x)^{2}+y^{2}}-
\sqrt{(1-x)^{2}+y^{2}}\right)
\end{aligned}
\qquad
\begin{aligned}
x &= uv \\
y &= \sqrt{(u^{2}-1)(1-v^{2})}
\end{aligned}
\EE
The equations for $ u,v $ are built from the distances to the points $ A, B $. They will only be considered in the upper half plane ($ \{(x,y) \mid y \geq 0\} $). By the equations \eqref{eq_ell_1} the upper half plane is homeomorphically mapped onto the region $ \{(u,v) \mid u \geq 1,\; -1 \leq v \leq 1\} $. The transformation back is recorded right in \eqref{eq_ell_1}. Outside the $ x $-axis the map $ (x,y) \mapsto (u,v) $ is diffeomorphic of class $ C^{\infty} $.

If the fundamental region in the first quadrant between 
$ \SS(M,1) $ and $ \SS(A,2) $ is transformed with the elliptic coordinates then the image lies in the square 
$ \{ (u,v) \mid 1 \leq u \leq 2, \; 0 \leq v \leq 1\} $ and is given there by
\BE
\label{eq_ell_2}
u^{2}+v^{2} \geq 2, \qquad u+v \leq 2.
\EE
The curves $ y = \psi_{1}(x) $, $ y = \psi_{2}(x) $ which separate the fundamental region into four parts are subsequently given by 
\BE
\label{eq_ell_3}
(1-uv)^2((u-v)^2+4(2-uv))-(u-v)^4 = 0, \qquad
(u^2-1)(1-v^2)-(1+uv)^2 = 0
\EE
with the corresponding assignment of the boundaries:

\begin{center}
\begin{minipage}{10.5cm}
\bild{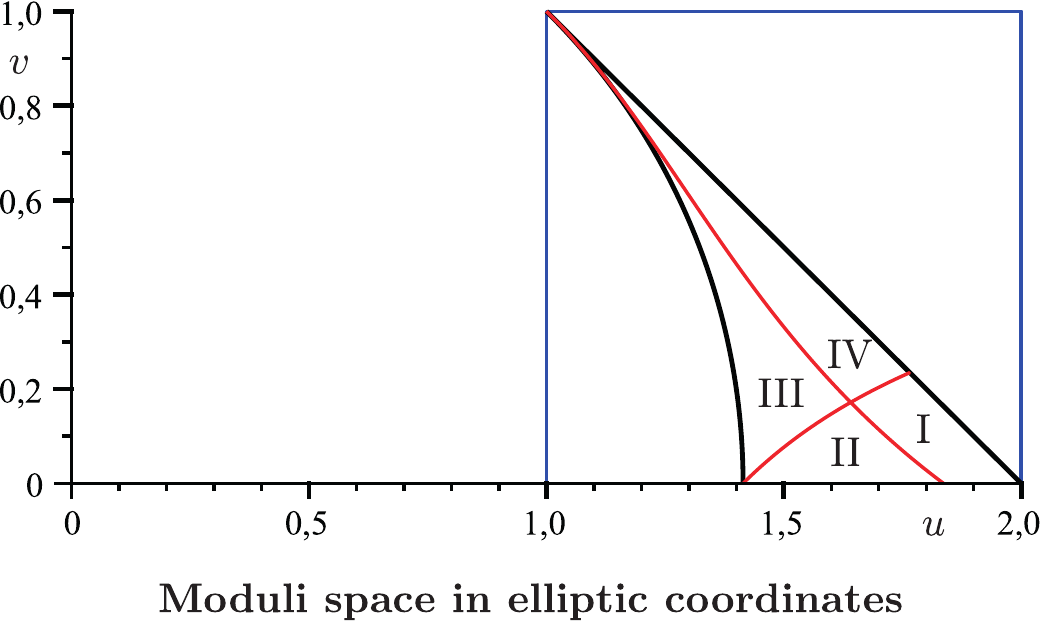}
\end{minipage}
\end{center}

The distinctive points on the $ u $-axis have the coordinates
\BE
\label{eq_ell_4}
\sqrt{2} = 1.4142\ldots, \quad
\sqrt{\frac{\sqrt{33}+1}{2}} = 1.8363\ldots 
\EE
The section point of the separating lines and the section point on the line $ u+v = 2 $ are
\BE
\label{eq_ell_5}
(1.6409\ldots,0.1716\ldots), \quad 
\left(1+\sqrt{2-\sqrt{2}},1-\sqrt{2-\sqrt{2}}\right) =
(1.7653\ldots,0.2346\ldots).
\EE
\TB{\TI{End of Intermezzo}}

All the regions will now be examined under the conditions
\BE
\label{eq_obb_21}
x > 0, \quad y > 0, \quad x^{2}+y^{2} > 1, 
\quad (1+x)^{2}+y^{2} < 4.
\EE
The regions I to IV are meant \TI{without} these points, however with inclusion of the separating curves $ y = \psi_{1}(x) $ and 
$ y = \psi_{2}(x) $. The comparisons of the $ d_{k} $ among themselves are done including the equality signs. The names I to IV will be maintained for the image regions in elliptic coordinates (using the corresponding relations of the boundaries).

\TI{For region} I:

Under the additional conditions
\BE
\label{eq_I_1}
y \geq \psi_{1}(x), \quad y \geq \psi_{3}(x)
\EE
one must find the minimum of the entries in the column `region I' of Table  \eqref{eq_obb_19}. 

The following picture shows a typical pattern for this region
\begin{center}
\begin{minipage}{7.5cm}
\bild{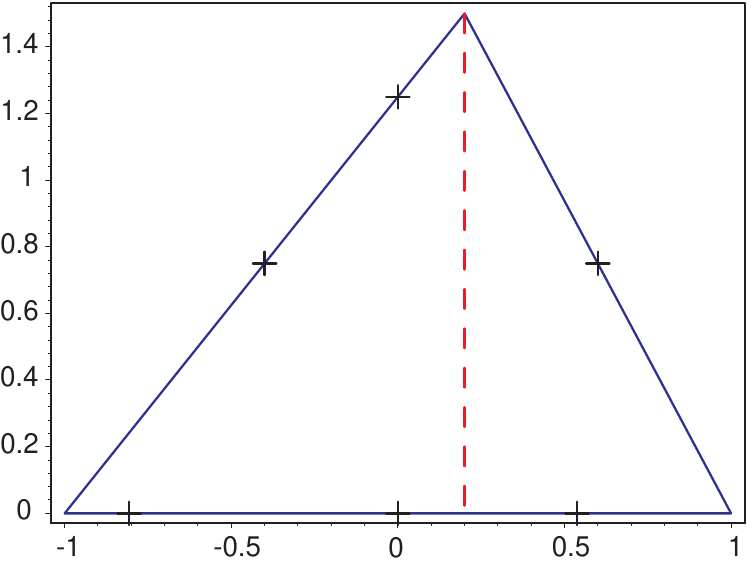}
\end{minipage}
\end{center}

Hereby, one may conjecture that the minimum is $ d_{3} = y $. This is really confirmed by the following comparisons:

\TI{Comparison with $ d_{1} $:} 
The inequality $ d_{3} \leq d_{1} $ reduces itself to 
\BE
\label{eq_I_2}
y \leq \sqrt{3}-\frac{x}{\sqrt{3}}.
\EE
This inequality describes the ordinate set of the tangent of  
$ \SS(A,2) $ at $ (0,\sqrt{3}) $, so is satisfied in all regions.

\TI{Comparison with $ d_{2} $:}
This is obvious (equality for $ x = 0 $).

\TI{Comparison with $ d_{4} $:}
The relation $ d_{3} \leq d_{4} $ is the inequality for the arithmetic and geometric mean for $ (1+x)^{2} $ and $ y^{2} $.

\TI{Comparison with $ d_{6} $:} 
This is obvious.

The strict inequality \eqref{eq_230808_65_1} says in the present case: 
\BE
\label{eq_I_3}
2+\sqrt{(1+x)^{2}+y^{2}}+\sqrt{(1-x)^{2}+y^{2}} > 2\sqrt{3}y
\EE
or, in the elliptic coordinates
\BE
\label{eq_I_4}
1+u > \sqrt{3}\;\sqrt{(u^{2}-1)(1-v^{2})}.
\EE
After squaring, this is equivalent to
\BE
\label{eq_I_5}
(u+1)(3uv^{2}-3v^{2}-2u+4) > 0.
\EE
By $ u+1 \geq 2 $ this is reduced to
\BE
\label{eq_I_6}
3(u-1)v^{2} > 2(u-2).
\EE
But this inequality is valid by observing $ v > 0 $ and 
$ 1 < u \leq 2 $.

\TI{For region} II:

Under the additional conditions
\BE
\label{eq_II_1}
\psi_{3}(x) \leq y \leq \psi_{1}(x)
\EE
one must find the minimum of the entries in the column `region II' of Table \eqref{eq_obb_19}. 

\begin{center}
\begin{minipage}{7.5cm}
\bild{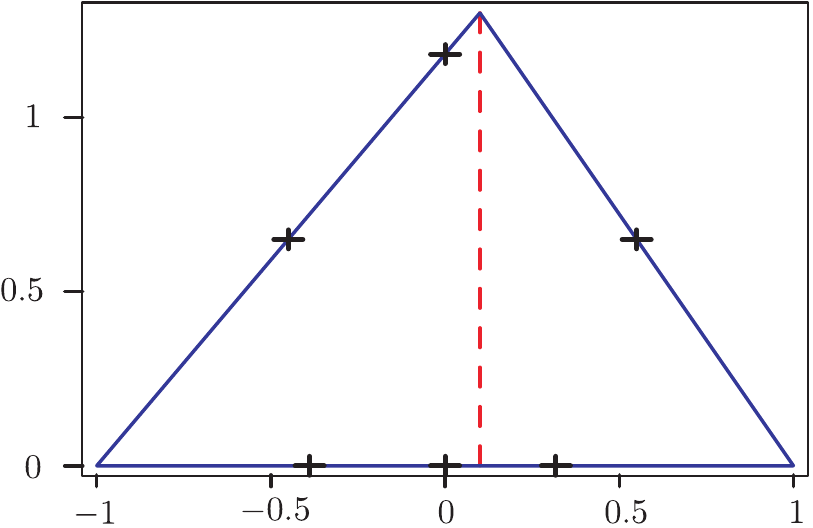}
\end{minipage}
\end{center}

This picture shows a typical pattern for this region, and again, one can conjecture that the minimum is $ d_{3} = y $. This is confirmed by the following comparisons.

\TI{Comparison with $ d_{1} $:} 
The inequality $ d_{3} \leq d_{1} $ reduces itself to the arithmetic-geometric inequality for $ (1-x)^{2} $ and $ y^{2} $. 

\TI{Comparison with $ d_{2} $, $ d_{4} $, $ d_{6} $:} 
This runs as for region I.

The strict inequality \eqref{eq_230808_65_1} is here expressed in the same way as for the region I, i.e. by \eqref{eq_I_3} and, equivalently, by  \eqref{eq_I_6}. Here, one has 
$ 2(u-2) < 2(1.9-2) = -0.2 $ by \eqref{eq_ell_4}, hence 
$ 3(u-1)v^{2} >  2(u-2) $.

\TI{For region} III:

Under the additional conditions
$$ 
y \leq \psi_{1}(x), \qquad y \leq \psi_{3}(x) 
$$
one must find the minimum of the entries in the column 
`region III' of Table \eqref{eq_obb_19}.

The following picture shows a typical pattern for this region:

\begin{center}
\begin{minipage}{7.7cm}
\bild{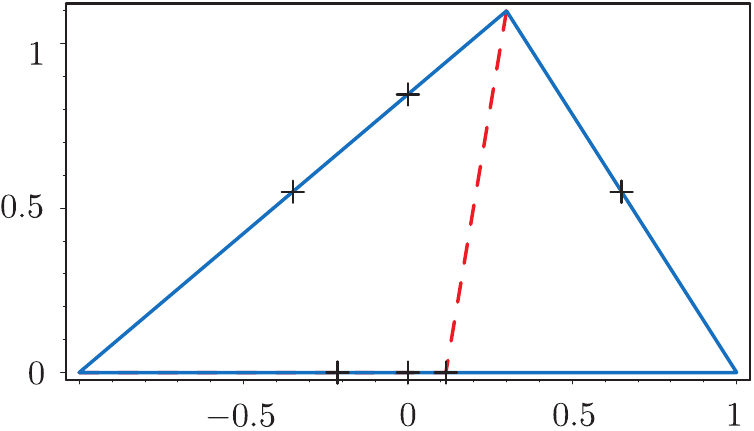}
\end{minipage}
\end{center}

Again, one can conjecture that the minimum is $ d_{4} $. This is confirmed by the following comparisons:

\TI{Comparison with $ d_{1} $:} 

The inequality $ d_{4} \leq d_{1} $ is equivalent to
\BE
\label{eq_III_2}
2x(x^{2}+y^{2}-1) \geq 0.
\EE
Obviously, this is true.

\TI{Comparison with $ d_{2} $:}

The inequality $ d_{4} \leq d_{2} $ is equivalent to
\BE
\label{eq_III_3}
(x^{2}+y^{2}-1)(3x^{2}+4x+1-y^{2}) \geq 0.
\EE
The first factor is positive, the second factor is $ \geq 0 $ if 
\BE
\label{eq_III_4}
y \leq \sqrt{(1+x)(3x+1)}.
\EE
By $ y \leq \psi_{3}(x) = 1+x $ this is certainly true if 
$$ 
1+x \leq \sqrt{(1+x)(3x+1)}.
$$
After squaring and cancelling of $ 1+x > 0 $ this is equivalent to 
$ x \geq 0 $. So, \eqref{eq_III_4} and then \eqref{eq_III_3} are 
satisfied.

The assertion \eqref{eq_230808_65_1} now says in the strict case:
\BE
\label{eq_III_5}
2+\sqrt{(1+x)^{2}+y^{2}}+\sqrt{(1-x)^{2}+y^{2}} > 
\sqrt{3}\;\frac{(1+x)^{2}+y^{2}}{1+x}
\EE
or, equivalently with elliptic coordinates
\BE
\label{eq_III_6}
1+u > \frac{\sqrt{3}}{2}\;\frac{(u+v)^{2}}{1+uv}.
\EE 
Since, in the region III, the $ y $-coordinate is always 
$ \leq \sqrt{2} $, it is sufficient for the proof of \eqref{eq_III_6} to assume $ y \leq y_{0} := \sqrt{2} $, what means in elliptic coordinates $ (u^{2}-1)(1-v^{2}) \leq 2 $  or
\BE
\label{eq_III_7}
v^{2} \geq 1-\frac{2}{u^{2}-1}{u^{2}-1}.
\EE
As long as the right hand side is $ \leq 0 $ (including 
$ = -\infty $) this is trivially true. This right hand side is strictly monotonously increasing as a function of $ u $ and has the unique zero $ \sqrt{3} $. So, one must show: From 
\BE
\label{eq_III_8}
u \geq \sqrt{3}, \qquad 
v \geq \sqrt{1-\frac{2}{u^{2}-1}} =: g(u)
\EE
follows, in the region III: \eqref{eq_III_6}, i.e.
\BE
\label{eq_III_9}
(1+u)(1+uv)-\frac{\sqrt{3}}{2}(u+v)^{2} > 0.
\EE
This inequality can explicitly solved for $ v $:
\BE
\label{eq_III_10}
f(u) < v < F(u)
\EE
with
\BE
\label{eq_III_11}
\begin{aligned}
f(u) &:= 
\frac{u}{\sqrt{3}}(u-(\sqrt{3}-1))-
\frac{u+1}{\sqrt{3}}\;\sqrt{(u-\sqrt{3})^{2}+2\sqrt{3}-3}
\\[1ex]
F(u) &:= 
\frac{u}{\sqrt{3}}(u-(\sqrt{3}-1))+
\frac{u+1}{\sqrt{3}}\;\sqrt{(u-\sqrt{3})^{2}+2\sqrt{3}-3}.
\end{aligned}
\EE
The function $ F $ is, for $ u \in [1,2] $, always $ > 1 $. For this, first consider the equation $ F(u) = 1 $ for all 
$ u \in \R $. The solutions of this equation result by equating 
the left hand side of \eqref{eq_III_9} to zero if $ v = 1 $, so from $ (u+1)^{2}(2-\sqrt{3}) = 0 $, and this is only true for 
$ u = -1 $. By $ F(2) = 1+\sqrt{3} $ one has $ F(u) > 1 $ for 
$ u \in [1,2] $. So, the second inequality in \eqref{eq_III_10} is automatically satisfied for $ 1 \leq u \leq 2 $, 
$ 0 \leq v \leq 1 $.

In $ [1,2] $, the function $ f $ changes its sign exactly twice. For this, consider the equation $ f(u) = 0 $ first for all 
$ u \in \R $. The solutions of this equation result by equating 
the left hand side of \eqref{eq_III_9} to zero if $ v = 0 $, so from $ \sqrt{3}u^{2}-2u-2 = 0 $. The solutions are 
\BE
\label{eq_III_12}
u_{0} := 
\frac{1}{\sqrt{3}}\left(1+2\sqrt{1+2\sqrt{3}}\right) = 1.7972\ldots,
\;
\frac{1}{\sqrt{3}}\left(1-2\sqrt{1+2\sqrt{3}}\right) =
-0.64249\ldots
\EE
By $ f(0) = -1 $, one has $ f(u) < 0 $ for 
$ u \in \left[1,u_{0}\right[ $ and $ f(u) > 0 $ for 
$ u \in \left]u_{0},2\right] $.

The graph of the function $ g $ from \eqref{eq_III_8} cuts the line $ u+v = 2 $ at
\BE
\label{eq_III_13}
u = u_{1} := 1+\sqrt{2-\sqrt{2}} = 1.7653\ldots
\EE
Thus, the points $ (u,v) $ with \eqref{eq_III_8} only can belong to the transformed moduli space if $ u \in [\sqrt{3},u_{1}] $. 
But in this interval, one has $ f(u) < 0 $ (since 
$ u_{1} < u_{0} $) so \eqref{eq_III_10} is satisfied there.
Hereby it is evident for the region III that the asserted inequality is valid in a strict form.

The situation is elucidated by the following picture:

\begin{center}
\begin{minipage}{16.3cm}
\bild{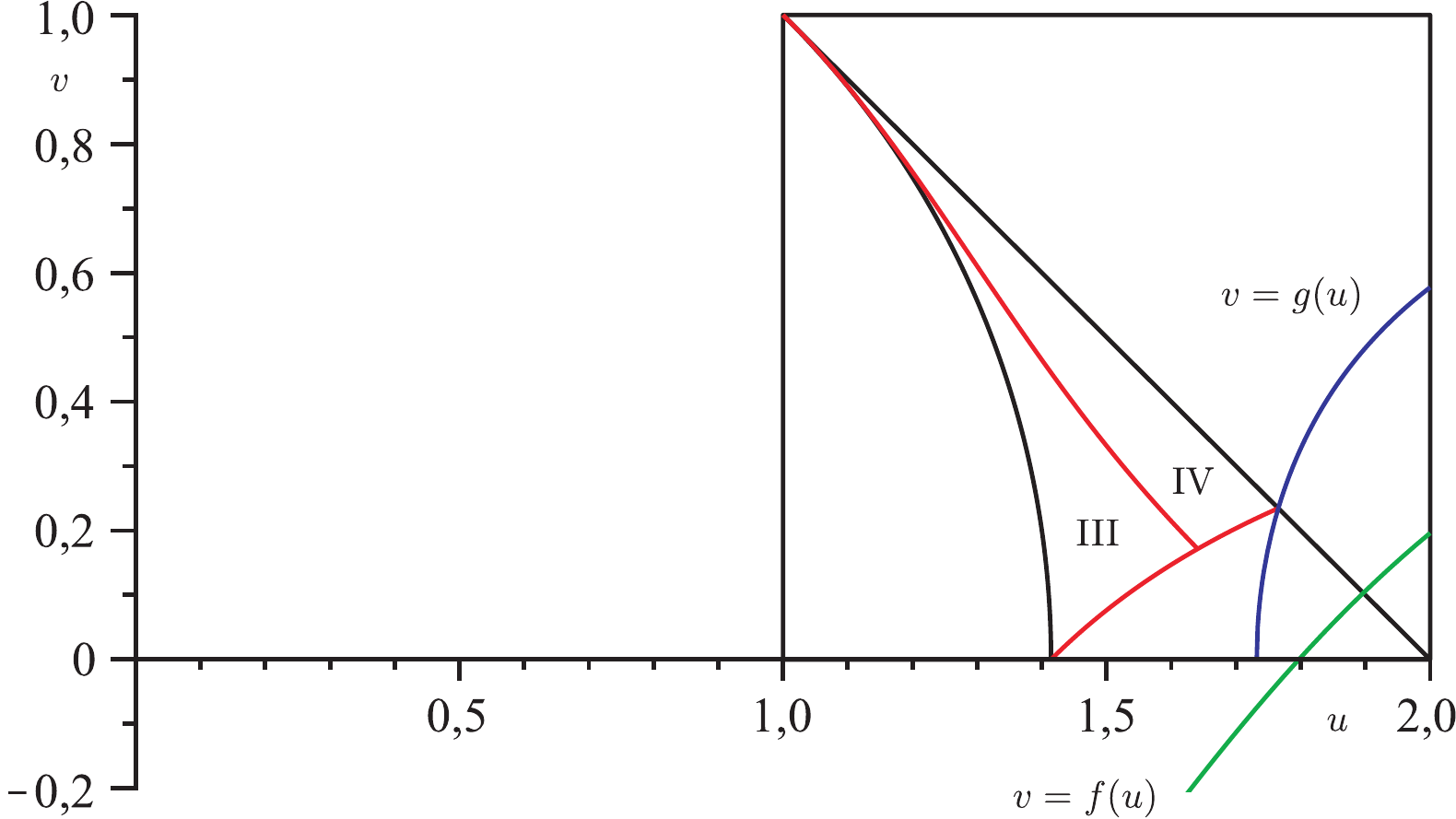}
\end{minipage}
\end{center}

and again under magnification around the points 
$ \sqrt{3}, u_{1}, u_{0} $ on the $ u $-axis:

\begin{center}
\begin{minipage}{9cm}
\bild{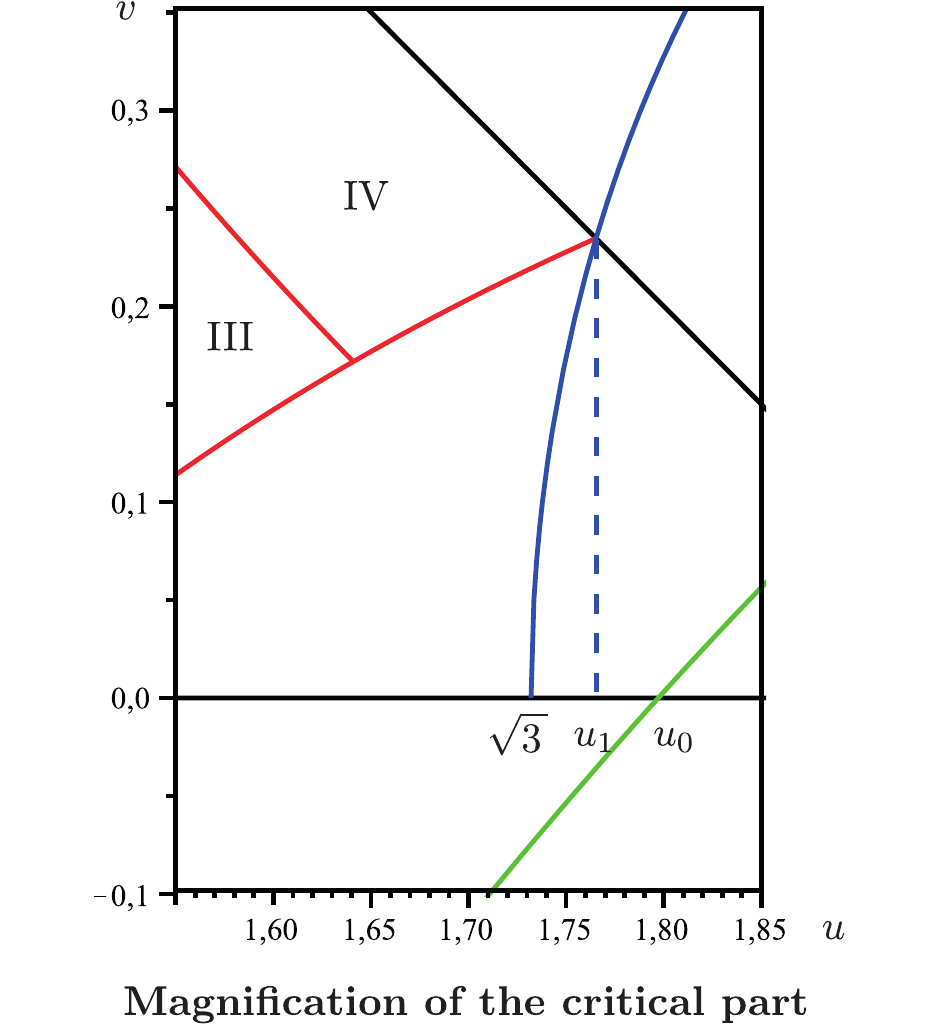}
\end{minipage}
\end{center}

In a moment, it will come out that the same picture is also appropriate for the last region IV.

\TI{For region} IV:

Under the additional conditions 
$$ 
\psi_{1}(x) \leq y \leq \psi_{3}(x) 
$$
one must find the minimum of the entries in the column 
`region IV' of Table \eqref{eq_obb_19}.

The following picture shows a typical pattern for this region.

\begin{center}
\begin{minipage}{7.7cm}
\bild{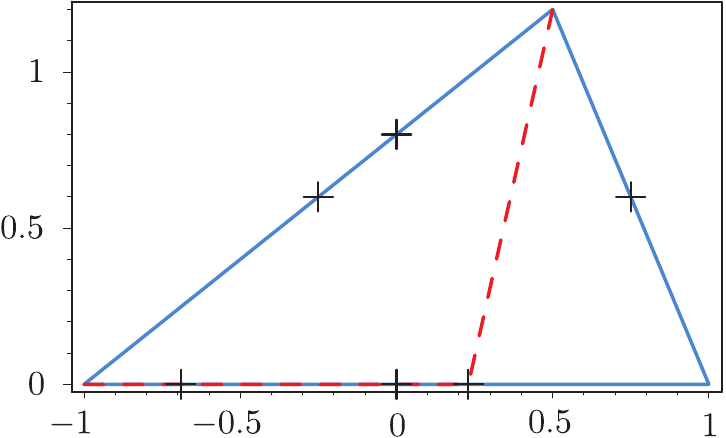}
\end{minipage}
\end{center}

One may conjecture that again $ d_{4} $ yields the minimum. This will be confirmed by the following comparisons:

\TI{Comparison with $ d_{1} $:}
The inequality $ d_{4} \leq d_{1} $ is equivalent to 
\BE
\label{eq_IV_1}
((1+x)^{2}+y^{2})^{2}-(1+x)^{2}((x-3)^{2}+y^{2}) \leq 0.
\EE
This is equivalent to
\BE
\label{eq_IV_2}
y^{2} \leq
\frac{1+x}{2}\left(\sqrt{x^{2}-30x+33}-(1+x)\right).
\EE
The expression in the large parentheses on the right turns out to be $ \geq 0 $ such that one has to confirm that in the region IV the following holds:
\BE
\label{eq_IV_3}
y \leq \sqrt{\frac{1+x}{2}}\;\sqrt{\sqrt{x^{2}-30x+33}-(1+x)}
=:h(x).
\EE
The graph of $ h $ over $ [0,1] $ mostly runs outside the moduli space. One has 
\BE
\label{eq_IV_4}
h(0) = \sqrt{\frac{\sqrt{33}-1}{2}} < \sqrt{3},
\EE
so $ (0,h(0)) \in \mathcal{M} $, but this graph cuts $ \SS(A,2) $ for $ 0 < x < 1 $ exactly once namely for
\BE
\label{eq_IV_5}
x = x_{0} := \frac{\sqrt{17}-3}{4} = 0.2807\ldots
\EE
This results by solving the equation $ h(x) = \sqrt{4-(1+x)^{2}} $ which is equivalent to $ (1-x)(2x^{2}+3x-1) = 0 $. So, \eqref{eq_IV_3} is satisfied for all $ x \in \mathcal{M} $ with 
$ x \geq x_{0} $. Since all points $ (x,y) $ in the closure of the region IV have an abscissa $ x \geq 0.2817\ldots $ it is evident that Eqn. \eqref{eq_IV_1} is valid there always. 

\TI{Comparison with $ d_{2} $:} 
This has already been confirmed with region III.

Further comparisons are not necessary, cf. Eqn. \eqref{eq_obb_19}.

As for region III it remains to prove the inequality \eqref{eq_III_5}. However, this runs exactly by the same arguments as for region III because the bound $ y_{0} = \sqrt{2} $ there has been chosen in such a way that all points in the closure of IV also have an ordinate $ y \leq y_{0} $. 

Finally, the cases excluded with \eqref{eq_obb_21} must be examined. For this one still has $ x^{2}+y^{2} > 1 $.

\TI{For $ x = 0 $, where $ 1 < y \leq \sqrt{3} $:}

Such a point $ (x,y) $ can be obtained by a limit process 
in the sense of Lemma \ref{lem_3.2}, part (i), of points 
in the interior of the region I or II, with the effect that 
$ \delta = y $ holds true:

\vspace{-10ex}
\begin{center}
\begin{minipage}{10.5cm}
\bild{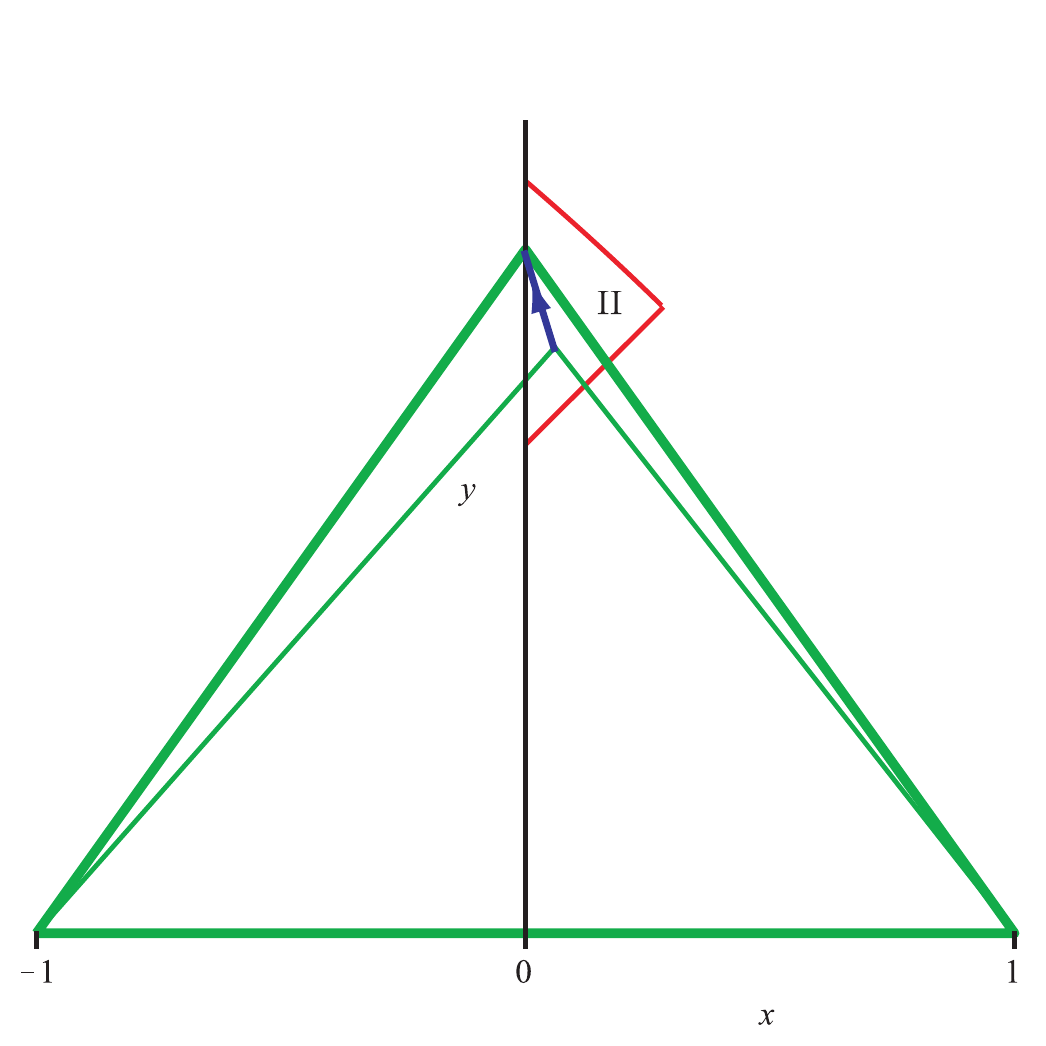}
\end{minipage}
\end{center}

In both cases one must show
\BE
\label{eq_L_1}
2+2\sqrt{1+y^{2}} \geq 2\sqrt{3}y
\EE
with equality exactly for $ y = \sqrt{3} $.

Equivalent with this is, by rearranging, 
$ 0 \geq 2y(y-\sqrt{3}) $. This is obvious with equality exactly for $ y = \sqrt{3} $.

\TI{For $ (x,y) \in \SS(A,2) $, where $ 0 \leq x < \sqrt{2}-1 $:}

Again, such a point $ (x,y) $ can be obtained as a limit of points in the interior of the region I, using Lemma \ref{lem_3.2}, part (i). Regarding $ \delta = y $ and 
$ (1+x)^{2}+y^{2} = 4 $, one has to show:
\BE
\label{eq_L_2}
2+2+\sqrt{(1-x)^{2}+4-(1+x)^{2}} \geq 2\sqrt{3}y
\EE
with equality exactly for $ (x,y) = (0,\sqrt{3}) $. Rearranging  \eqref{eq_L_2} yields 
\BE
\label{eq_L_3}
\frac{1}{16}x(x+3)(8-3x-9x^{2}) \geq 0.
\EE
The zeros of the last polynomial factor are
\BE
\label{eq_L_4}
\frac{-1-\sqrt{33}}{6} = -1.1240\ldots, \qquad
\frac{-1+\sqrt{33}}{6} = 0.7907\ldots.
\EE
So, this factor is positive for 
$ 0 \leq x < \sqrt{2}-1 = 0.4142\ldots $, and Eqn. \eqref{eq_L_3} is fulfilled with equality exactly for 
$ x = 0 $ and thus $ y = \sqrt{3} $.

\TI{For $ (x,y) \in \SS(A,2) $, where $ \sqrt{2}-1 \leq x < 1 $:}

Such a point $ (x,y) $ can be obtained as a limit of points in the interior of the region IV such that $ \delta = d_{4} $ holds, again using Lemma \ref{lem_3.2}, part (i). By 
$ (1+x)^{2}+y^{2} = 4 $ one has to show
\BE
\label{eq_L_5}
2+2+\sqrt{4-4x} > \sqrt{3}\;\frac{4}{1+x},
\EE
equivalently
\BE
\label{eq_L_6}
\sqrt{1-x} > \frac{2(\sqrt{3}-1)-2x}{1+x}.
\EE
Here, the right hand side is negative for $ x > \sqrt{3}-1 $ such that Eqn. \eqref{eq_L_6} has only to be confirmed for 
$ \sqrt{2}-1 \leq x \leq \sqrt{3}-1 $. In this case, the assertion 
\eqref{eq_L_6} turns out to be equivalent to 
$ 0 > x^{2}-2x+2\sqrt{3}-3 $. The last polynomial has the zeros
\BE
\label{eq_L_7}
2-\sqrt{3} = 0.2679\ldots, \qquad \sqrt{3} = 1.7320\ldots,
\EE
so is negative in between, in particular for 
$ x \in [\sqrt{2}-1,\sqrt{3}-1] = [0.3142\ldots,0.7320\ldots] $.

So everything has been shown.
\QED

\medskip
\TB{\TI{The case of quadrangles}}

Here, we only can offer a strong conjecture which causes a certain surprise because the extreme figure is definitely not the square.

\BX
\label{ex_unviex}
For a \TI{square} $ P $ one has
\BE
\label{eq_unviex_i}
\frac{L(P)}{\delta(P)} = \frac{8}{5}\sqrt{5} = 
3.577\;708\;763\ldots.
\EE
See the following picture.
\EX

For a \TI{general quadrangle} $ P $ we may choose the line from 
$ B = (-1,0) $ to $ E = (1,0) $ as a diameter $ D $ (see Lemma 
\ref{lem_diam}). 

If $ D $ is not an edge of $ P $ we may place 
one more vertex as $ A = (u0,u) $ with $ -1 < u_0 < 1 $, $ u > 0 $ in the upper half plane and the last vertex $ C = (v0,-v) $ with $ -1 < v_0 < 1 $, $ v > 0 $ in the lower half plane. Then a simple search programme shows an approximate minimum of the quotient 
$ L(P)/\delta(P) $ at 
$$ 
u_0 = v_0 = 0,\quad u := 1.472, \quad v = 0.217
$$
So, likely, the minimum figure is symmetric w.r.t. the 
$ y $-axis. Moreover, the algorithm from above shows that the 
distinguished chords all have the same length (and are grounding orthogonally on the receiving edges).

Now, an exact search can be started, considering only the 
fixed vertices $ A $ and $ E $ and the variable vertices 
$ A = (0,u) $ and $ C = (0,-v) $, thus reducing the problem from dimension $ 4 $ to dimension $ 2 $, with 
$ (u,v) \in \Rplus\cross\Rplus $ and the side condition 
$ u \neq v $ ($ u = v $ would yield a square). Also, the diameter of length $ 2 $ between $ A $ and $ E $ results in some additional restrictions.

The equality of the two perpendiculars from $ A $ to $ B \vee C $, resp. from $ B $ to $ A \vee E $ is expressed by the relation
$$
\frac{1}{2}\,\frac{u+v}{\sqrt{v^2+1}} = \frac{u}{\sqrt{u^2+1}}
$$
which comes down algebraically to
$$
(u^3-3u+(3u^2-1)v)(u-v) = 0.
$$
or, using $ u-v \neq 0 $, to 
\BE
\label{vu}
v = \frac{3-u^2}{3u^2-1}\,u.
\EE
where $ 0 < u < 1 $ and $ v > 0 $ require $ u > 1/\sqrt{3} $.

So the problem is reduced to the minimum question of a real function, namely of the quotient 
`perimeter/length of the perpendicular from $ B $ to 
$ A \vee E $', i.e.
\BE
\label{vu1}
\frac{2\sqrt{u^2+1}+2\sqrt{v^2+1}}{\ds \frac{u}{\sqrt{u^2+1}}},
\EE
with $ v $ inserted according to  \eqref{vu}. This function has the representation
$$
f(u) := 4u\,\frac{u^2+1}{3u^2-1}, \quad 
\frac{1}{\sqrt{3}} < u < 1.
$$
Its minimum condition is $ 3u^4-6u^2-1 = 0 $, with the only positive solution
$$
u = \frac{\sqrt{3}}{3}\,\sqrt{3+2\sqrt{3}}.
$$
From this one derives the corresponding $ v $ with
\eqref{vu} and the isoperimetric quotient with \eqref{vu1}, and one can state:

If $ P $ is the \TB{kite} with the vertices
\BE
\label{eq_150808_15_1}
(-1,0), \quad (1,0), 
\quad \left(0,\frac{\sqrt{3}}{3}\,\sqrt{2\sqrt{3}+3}\right),  
\quad
\left(0,-\frac{1}{3}\,\sqrt{2\sqrt{3}-3}\right),
\EE
then 
\BE
\label{eq_150808_15_2}
\frac{L(P)}{\delta(P)} = 
\frac{4}{3}\,\sqrt{2\sqrt{3}+3} = 3.389\;946\;342\ldots
\EE
Due to these weird number relations, this figure $ P $ shall be called a \TB{magic kite}. 

\vspace{-2ex}
\hspace*{-9mm}
\begin{minipage}{16cm}
\bild{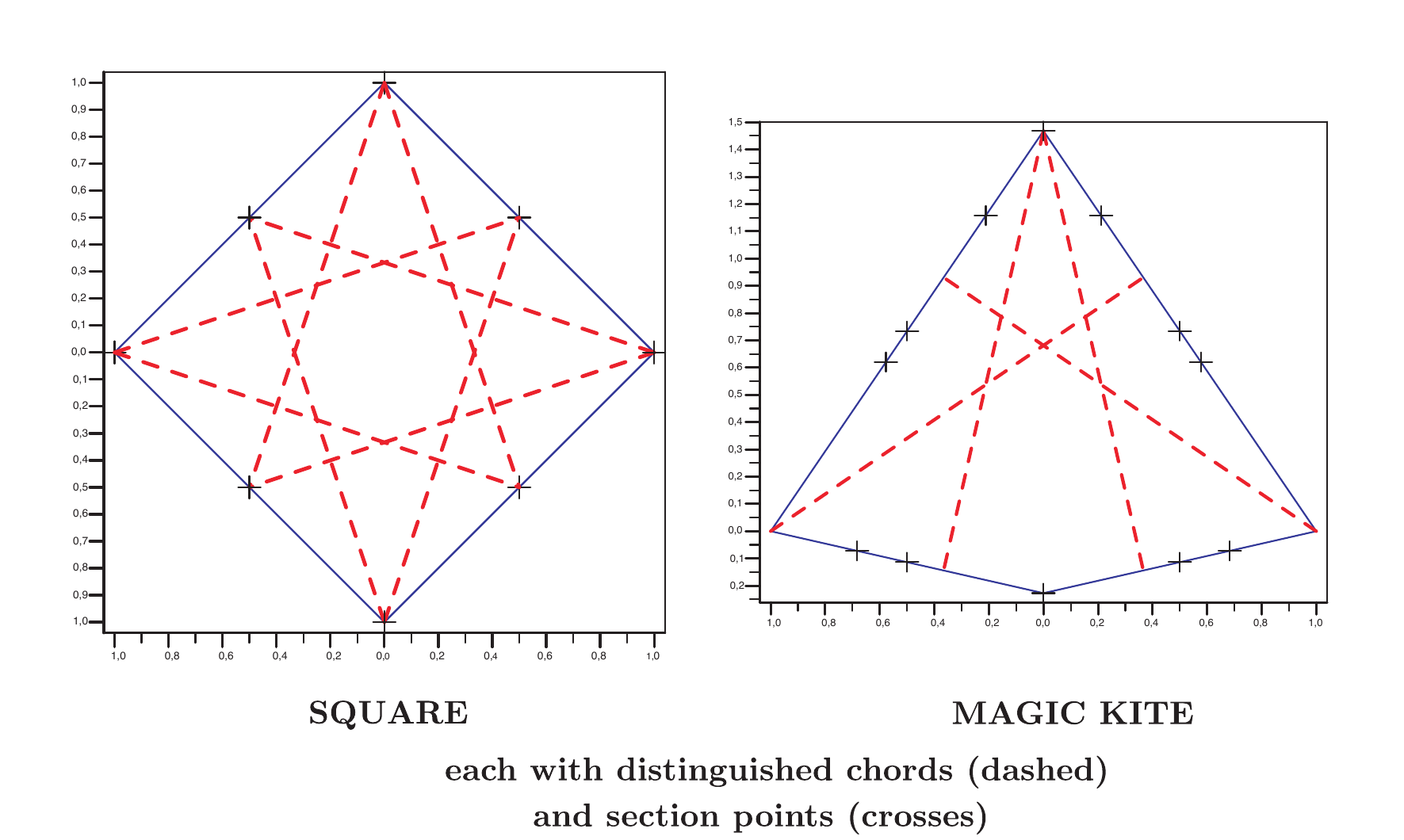}
\end{minipage}

If $ D $ is an edge, the two other vertices $ A $ and $ C $ of $ P $ lie in the same half plane, say in the upper one. Numerical optimization programmes don't give any indication that, varying 
$ A $ and $ C $, could yield a genuine convex quadrangle with minimal invariant $ \delta $.

These arguments provide a strong support for the

\TB{Conjecture 2:}
\TI{For any convex quadrangle $ P \subset \R^2 $ there holds the inequality} 
$$
L(P) \geq \frac{4}{3}\,\sqrt{2\sqrt{3}+3}\cdot\delta(P)
$$
\TI{between the perimeter $ L(P) $ and the invariant 
$ \delta(P) $, with equality exactly for magic kites.}

\bigskip 
\TB{\large References}

\TI{Amir, D., and Ziegler, Z. [1980]:}
Relative Chebyshev centers in normed linear spaces I. 
J. Approximation Theory \TB{29} (1980), 235-252

\TI{Bangert, V. [1981]:}
Totally convex sets in complete Riemannian manifolds. 
J. Diff. geom. \TB{16} (1981), 333-345

\TI{Bonnesen, T., und Fenchel, W. [1934]:}
Theorie der konvexen K\"orper. 
I-VII u. 1-164. 
Springer Berlin \TB{1934}. 
(Reprint Chelsea Publ. Comp. New York 1971)

\TI{Kobayashi, S. and Nomizu, K. [1963]:} 
Foundations of differential geometry, Vol. I. 
i-xi a.~1-329. 
Interscience Publishers John Wiley and Sons \TB{1963}

\TI{Phelps, R.R. [1957]:} 
Convex sets and nearest points. 
Proc. Amer. Math. Soc. \TB{8} (1957), 790-797

\TI{Wortmann, Ch. [2004]:} 
Algorithmische Analyse von rohrf\"ormigen Fl\"achen f\"ur das 
Reverse Engineering. 
Thesis preprint Fachbereich Informatik Dortmund \TB{2004}

\bigskip
\TB{Author's address:} \\
Fakult\"at f\"ur Mathematik \\
Technische Universit\"at Dortmund\\
D-44221 Dortmund, Germany \\
e-mail: rolf.walter@tu-dortmund.de

\end{document}